\documentclass[a4paper,10pt,twoside,bibliography=totoc]{scrbook}

\title{Representation of Conditional Expectations in Gaussian Analysis on Sequence Spaces}
\author{Felix Riemann}

\usepackage{hyperref}

\usepackage{ucs}
\usepackage[utf8x]{inputenc}

\usepackage{geometry}
\usepackage[T1]{fontenc}

\usepackage[english]{babel}

\usepackage{amsmath}
\usepackage{amsthm}
\usepackage{amssymb}
\usepackage{amsfonts}

\usepackage{graphicx}

\theoremstyle{plain}
\newtheorem{theorem}{Theorem}[chapter]
\newtheorem{lemma}[theorem]{Lemma}
\newtheorem{proposition}[theorem]{Proposition}
\newtheorem{corollary}[theorem]{Corollary}

\theoremstyle{definition}
\newtheorem{definition}[theorem]{Definition}
\newtheorem{example}[theorem]{Example}
\newtheorem{application}[theorem]{Application}

\theoremstyle{remark}
\newtheorem{remark}[theorem]{Remark}
\newtheorem*{notation}{Notation}

\newcommand{\E}{\mathbb{E}}
\newcommand{\K}{\mathbb{K}}
\newcommand{\C}{\mathbb{C}}
\newcommand{\R}{\mathbb{R}}
\newcommand{\N}{\mathbb{N}}
\newcommand{\Z}{\mathbb{Z}}

\newcommand{\calF}{\mathcal{F}}
\newcommand{\calG}{\mathcal{G}}
\newcommand{\calM}{\mathcal{M}}
\newcommand{\calN}{\mathcal{N}}
\newcommand{\calH}{\mathcal{H}}
\newcommand{\calP}{\mathcal{P}}
\newcommand{\calW}{\mathcal{W}}

\newcommand{\wh}{\widehat}
\newcommand{\wot}{{\wh{\otimes}}}
\newcommand{\otn}{{\otimes n}}
\newcommand{\otm}{{\otimes m}}
\newcommand{\wotn}{{\wot n}}
\newcommand{\wotm}{{\wot m}}
\newcommand{\wickpoly}[2][n]{\big\la #2,\colon\cdot^{\otimes #1}\colon\big\ra}
\newcommand{\wickpolyo}[2][n]{\big\la #2,\colon\omega^{\otimes #1}\colon\big\ra}
\newcommand{\indi}{\mathbf{1}}
\newcommand{\la}{\langle}
\newcommand{\ra}{\rangle}

\renewcommand{\phi}{\varphi}
\renewcommand{\epsilon}{\varepsilon}

\DeclareMathOperator{\Id}{Id}
\DeclareMathOperator{\id}{id}
\DeclareMathOperator{\dist}{dist}
\DeclareMathOperator{\dx}{dx}
\DeclareMathOperator{\spann}{span}

\begin{document}

\begin{titlepage}
\begin{center}

\vspace*{30mm}
\textsc{\LARGE -- Diploma Thesis --}

\vspace{20mm}
\hrule height 2pt
\vspace{5mm}
{\huge\bfseries \parbox[0pt][1.3em][c]{0cm}{}Representation of Conditional \parbox[0pt][1.3em][c]{0cm}{}Expectations in Gaussian Analysis on \parbox[0pt][1.3em][c]{0cm}{}Sequence Spaces}
\vspace{5mm}
\hrule height 2pt

\vspace{15mm}
{\large
 Handed in: March 31, 2011\\[2mm]
 Last modification: \today
}
\vspace{10mm}

\begin{minipage}{0.3\textwidth}
\begin{flushleft}
\large
\emph{Author:}\\
\textsc{Felix Riemann}
\end{flushleft}
\end{minipage}
\begin{minipage}{0.493\textwidth}
\begin{flushright}
\large
\emph{Advisor:}\\
\textsc{Prof. Dr. Martin Grothaus}
\end{flushright}
\end{minipage}

\vspace{5mm}
\begin{minipage}{0.8\textwidth}
{\large\textsc{TU Kaiserslautern, Department of Mathematics}}\linebreak
{\large\textsc{Functional Analysis and Stochastic Analysis Group}}
\end{minipage}

\vspace{20mm}
\begin{figure}[h]
\centering
\includegraphics[width=0.4\textwidth]{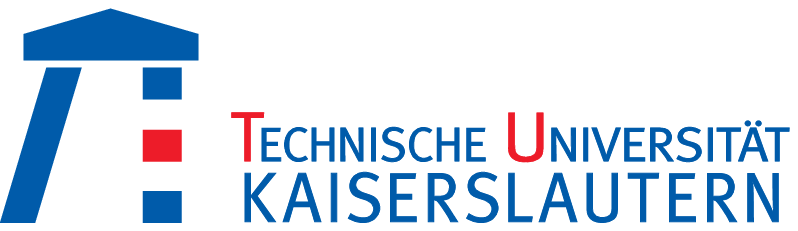}
\end{figure}

\end{center}
\end{titlepage}

\tableofcontents

\chapter*{Introduction}
\addcontentsline{toc}{chapter}{Introduction}

In infinite-dimensional analysis, the concept of Gaussian analysis is used to establish a Gaussian measure on a linear space of infinite dimension. While there does not exist a Lebesgue measure (i.e. a translation-invariant and locally finite measure which is non-trivial) in infinite dimensions, Gaussian measures can still be defined in this setting. Since Gaussian systems often occur in stochastic analysis, this theory becomes useful in solving stochastic partial differential equations. Furthermore, the area of White Noise Analysis may be used to approach Feynman path integrals in quantum physics.
\medskip

For a complete nuclear space $\calN$ which is densely and continuously embedded in a separable real Hilbert space $\big(\calH,(\cdot,\cdot)_\calH\big)$, a Gaussian measure can be constructed on $\calN'$, the topological dual space of $\calN$. By identifying $\calH$ with its topological dual space $\calH'$, we obtain a so-called Gel'fand triple $\calN\subset\calH=\calH'\subset\calN'$. The Bochner-Minlos theorem states that if we equip $\calN'$ with its cylindrical $\sigma$-algebra, for a characterisic function $C:\calN\to\C$ there exists a unique measure $\mu$ on $\calN'$ such that for all $\eta\in\calN$ one has
\[ \int_{\calN'}\exp(i\la\eta,\omega\ra)d\mu(\omega)=C(\eta), \]
where $\la\cdot,\cdot\ra$ denotes the canonical dual pairing between $\calN$ and $\calN'$. The characteristic function $C$ is usually chosen to be $\eta\mapsto\exp\left(-\frac{1}{2}(\eta,\eta)_\calH\right)$ in order to obtain a standard Gaussian measure. However, our goal is to insert a covariance operator $A$ in this measure.
\medskip

Motivated by applications (see Section~\ref{appsect}) we will build our Gaussian space around the sequence space $\ell^2(\calH)$. For this purpose we will characterize $\ell^2(\calH)$ and show that an operator $A\in L(\ell^2(\R))$ naturally can be extended to an element of $L(\ell_2(\calH))$. Such operators will be used as a special type of correlation operators which allow a more explicit representation of certain conditional expectations.
\medskip

In the second chapter the correlated Gaussian measure will be constructed. We will define $s(\calN)$, a dense nuclear subspace of $\ell^2(\calH)$, and establish the measure $\mu_A$ on its topological dual space $s'(\calN)$ by means of the Bochner-Minlos theorem, using the characterisic function
\[ s(\calN)\ni\phi\longmapsto\exp\left(-\frac{1}{2}(\phi,A\phi)_\calH\right)\in\R, \]
where the corvariance operator $A$ is a self-adjoint and positive definite operator in $L(\ell^2(\calH))$. We will rederive the orthogonal decomposition into wick polynomials in $L^2(\mu_A)$, called chaos decomposition, and see that it slightly differs from the usual decomposition, as the kernels are going to be elements from $\ell_A^2(\calH)$, the completion of $\ell^2(\calH)$ with respect to the inner product generated by $A$.
\medskip

Finally, in the third chapter, a representation for the conditional expectation of arbitrary random variables in $L^2(\mu_A)$ will be given, where we condition on a $\sigma$-algebra generated by monomials. Later on we will exploit the structure of the sequence space to obtain a representation of conditional expectations of monomials conditioned on countably many monomials and give an application of the results.
\medskip

Throughout the whole thesis, the Gaussian spaces will only be considered over the field $\R$, even if Gaussian analysis usually deals with the complexification of these spaces. However, the results in here may be generalized to cover the complex case. Furthermore, to simplify the proofs, we will consider the Hilbert space $\calH$ to satisfy $\dim\calH=\infty$, even though the proofs for finite-dimensional $\calH$ work similar.

\begin{chapter}{Square Summable Sequences}

Throughout all chapters, $\big(\calH,(\cdot,\cdot)_\calH\big)$ will be a real separable Hilbert space with $\dim \calH=\infty$. By $(e_k)_{k\in\N}$ we will denote the unit vectors in $\ell^2(\R)$, i.e. $e_k=(\delta_{k,l})_{l\in\N}\in\ell^2(\R)$ for $k\in\N$.

\begin{section}{Characterization}

\begin{definition}
 The Hilbert space of square summable sequences in $\calH$ we denote by
\[ \ell^2(\calH):=\left\{f\in \calH^\N:\sum_{k=1}^\infty\|f_k\|_\calH^2<\infty\right\} \]
together with its inner product
\[ \ell^2(\calH)\times\ell^2(\calH)\ni(f,g)\longmapsto(f,g)_{\ell^2(\calH)}:=\sum_{k=1}^\infty(f_k,g_k)_\calH\in\R. \]
The norm on $\ell^2(\calH)$ we denote by $\|\cdot\|_{\ell^2(\calH)}:=\sqrt{(\cdot,\cdot)_{\ell^2(\calH)}}$ as usual.
\end{definition}

\begin{lemma}
 The mappings
\[ \calH\times\ell^2(\R)\ni (h,x)\longmapsto h\bullet x:=(hx_k)_{k\in\N}\in\ell^2(\calH)\quad\text{and} \]
\[ \ell^2(\calH)\times\ell^2(\R)\ni (f,x)\longmapsto [f,x]:=\sum_{k=1}^\infty x_kf_k\in\calH \]
are bilinear and satisfy
\begin{enumerate}
 \item $\|\cdot\bullet x\|_{L(\calH,\ell^2(\calH))}=\|x\|_{\ell^2(\R)}$,
 \item $\|h\bullet\cdot\|_{L(\ell^2(\R),\ell^2(\calH))}=\|h\|_{\calH}$,
 \item $\|[\cdot,x]\|_{L(\ell^2(\calH),\calH)}=\|x\|_{\ell^2(\R)}$ and
 \item $\|[f,\cdot]\|_{L(\ell^2(\R),\calH)}\le\|f\|_{\ell^2(\calH)}$
\end{enumerate}
for all $x\in\ell^2(\R)$, $h\in\calH$ and $f\in\ell^2(\calH)$.
\begin{proof}
 The bilinearity of both mappings is clear. For $x\in\ell^2(\R)$ and $h\in\calH$ we have
\[ \|h\bullet x\|_{\ell^2(\calH)}^2=\sum_{k=1}^\infty x_k^2\|h\|_\calH^2=\|x\|_{\ell^2(\R)}^2\|h\|_\calH^2, \]
so $\cdot\bullet\cdot$ is well-defined and both \textit{(i)} and \textit{(ii)} are proven. For $f\in\ell^2(\calH)$ by the triangle inequality for $\|\cdot\|_\calH$ and the Cauchy-Bunyakovsky-Schwarz inequality we obtain
\begin{equation}\label{somelabel0}
 \left\|\sum_{k=m}^nx_kf_k\right\|_\calH\le\sum_{k=m}^n|x_k|\|f_k\|_\calH\le\left(\sum_{k=m}^n|x_k|^2\right)^\frac{1}{2}\left(\sum_{k=m}^n\|f_k\|_\calH^2\right)^\frac{1}{2}
\end{equation}
for all $n,m\in\N$. This yields that $\big(\sum_{k=1}^n x_kf_k\big)_{n\in\N}$ is a Cauchy-Sequence in $\calH$, so $[\cdot,\cdot]$ is well-defined. By setting $m=1$ in \eqref{somelabel0} and taking the limit $n\to\infty$ we obtain the estimation $\|[f,x]\|_\calH\le\|f\|_{\ell^2(\calH)}\|x\|_{\ell^2(\R)}$ which proves \textit{(iv)} and shows $\|[\cdot,x]\|_{L(\ell^2(\calH),\calH)}\le\|x\|_{\ell^2(\R)}$. For the proof of \textit{(iii)} it is left to show $\|[\cdot,x]\|_{L(\ell^2(\calH),\calH)}=\|x\|_{\ell^2(\R)}$. We compute
\[ \big\|[h\bullet x,x]\big\|_\calH=\left\|\sum_{k=1}^\infty hx_k^2\right\|_\calH=\left|\sum_{k=1}^\infty x_k^2\right|\|h\|_\calH=\|x\|_{\ell^2(\R)}^2\|h\|_\calH=\|x\|_{\ell^2(\R)}\|h\bullet x\|_{\ell^2(\calH)}. \]
This also shows that the inequality in \textit{(iv)} is sharp. If $f=(f_1,f_2,0,0,\dots)\in\ell^2(\calH)$ for some orthonormal $f_1,f_2\in\calH$, then $\|f\|_{\ell^2(\calH)}^2=2$ and for all $x\in\ell^2(\R)$ with $x\not=0$ we have
\[ \big\|[f,x]\big\|_\calH^2=\|x_1f_1+x_2f_2\|_\calH^2=x_1^2+x_2^2\le\|x\|_{\ell^2(\R)}^2<\|x\|_{\ell^2(\R)}^2\|f\|_{\ell^2(\calH)}^2, \]
so the inequality in \textit{(iv)} even may be strict.
\end{proof}
\end{lemma}

\begin{proposition}\label{rechenregeln}
 Let $f\in\ell^2(\calH)$, $h,g\in\calH$ and $x,y\in\ell^2(\R)$. We have the following identities:
\begin{enumerate}
 \item $(h\bullet x,g\bullet y)_{\ell^2(\calH)}=(h,g)_\calH(x,y)_{\ell^2(\R)}$.
 \item $(f,h\bullet x)_{\ell^2(\calH)}=\big([f,x],h\big)_\calH$.
 \item $[h\bullet x,y]=(x,y)_{\ell^2(\R)}h$.
\end{enumerate}
\begin{proof}
 This is done by the following straightforward computations:
\begin{enumerate}
 \item $(h\bullet x,g\bullet y)_{\ell^2(\calH)}=\sum_{k=1}^\infty(hx_k,gy_k)_\calH=\sum_{k=1}^\infty x_ky_k(h,g)_\calH=(h,g)_\calH(x,y)_{\ell^2(\R)}$.
 \item $(f,h\bullet x)_{\ell^2(\calH)}=\sum_{k=1}^\infty(f_k,hx_k)_\calH=\big(\sum_{k=1}^\infty x_kf_k,h\big)_\calH=\big([f,x],h\big)_\calH$.
 \item $[h\bullet x,y]=\sum_{k=1}^\infty hx_ky_k=(x,y)_{\ell^2(\R)}h$.
\end{enumerate}
\end{proof}
\end{proposition}

\begin{remark}
 For the sake of notational simplicity, we will often omit the indexes of the norms, i.e. we will write $\|\cdot\|$ for $\|\cdot\|_{\ell^2(\R)}$, $\|\cdot\|_\calH$, $\|\cdot\|_{\ell^2(\calH)}$ et cetera, since there is no risk of confusion. Analogously we will deal with inner products.
\end{remark}

\begin{proposition}\label{finseqdense}
 Let $(h_i)_{i\in\N}$ be an orthonormal basis of $\calH$. Then $\{h_i\bullet e_k:i,k\in\N\}$ is an orthonormal basis of $\ell^2(\calH)$.
\begin{proof}
 For $i,i',k,k'\in\N$ we clearly have
\[ (h_i\bullet e_k,h_{i'}\bullet e_{k'})=(h_i,h_{i'})(e_k,e_{k'})=\delta_{i,i'}\delta_{k,k'}=\delta_{(i,k),(i',k')} \]
by Proposition~\ref{rechenregeln}. Now let $f\in\ell^2(\calH)$ and $\epsilon>0$ be arbitrary. Choose $N\in\N$ such that $\sum_{k=N+1}^\infty\|f_k\|^2<\epsilon/2$. For each $k=1,\dots,N$ there exists $g_k\in\spann\{h_i:i\in\N\}$ with $\|f_k-g_k\|^2<\epsilon/2N$. Then for $g:=\sum_{k=1}^Ng_k\bullet e_k\in\spann\{h_i\bullet e_k:i,k\in\N\}$ we have
\[ \|f-g\|^2=\sum_{k=1}^\infty\|f_k-g_k\|^2=\sum_{k=1}^N\|f_k-g_k\|^2+\sum_{k=N+1}^\infty\|f_k\|^2<\epsilon. \]
\end{proof}
\end{proposition}

\begin{remark}
 While it is clear that for $f\in\ell^2(\calH)$ it holds $f=\sum_{k=1}^\infty f_k\bullet e_k$, it may not be too obvious that such an identity also exists for an arbitrary orthonormal basis $(b_k)_{k\in\N}$ of $\ell^2(\R)$, i.e. that there exists a sequence $(f_k^b)_{k\in\N}$ in $\calH$ such that $f=\sum_{k=1}^\infty f_k^b\bullet b_k$. The rest of this section will deal with the proof that this identity is valid for $f_k^b:=[f,b_k]$.
\end{remark}

\begin{definition}
 For an orthonormal basis $b=(b_k)_{k\in\N}$ in $\ell^2(\R)$ and $N\in\N$ we define
\[ \ell^2(\calH)_b^{(N)}:=\spann\left\{h\bullet b_k:h\in\calH,k=1,\dots,N\right\}\subset\ell^2(\calH). \]
\end{definition}

In the following, $b=(b_k)_{k\in\N}$ will always be an arbitrary orthonormal basis of $\ell^2(\R)$, if not stated otherwise.

\begin{proposition}
 For $N\in\N$ we have 
\[ \ell^2(\calH)_b^{(N)}=\left\{\sum_{k=1}^Nh_k\bullet b_k:h_1,\dots,h_N\in\calH\right\}. \]
\begin{proof}
 If $f\in\ell^2(\calH)_b^{(N)}$ then there exist $m\in\N$ and $\alpha_i\in\R$, $g_i\in\calH$ and $k_i\in\{1,\dots,N\}$ for $i=1,\dots,m$ such that
\[ f=\sum_{i=1}^m\alpha_ig_i\bullet b_{k_i}=\sum_{k=1}^N\Bigg(\sum_{\substack{i=1,\dots,m:\\k_i=k}}\alpha_ig_i\Bigg)\bullet b_k=\sum_{k=1}^Nh_k\bullet b_k,\quad\text{where }h_k:=\sum_{\substack{i=1,\dots,m:\\k_i=k}}\alpha_ig_i\in\calH. \]
\end{proof}
\end{proposition}

\begin{corollary}\label{finseqclosed}
 The space $\ell^2(\calH)_b^{(N)}$ is a closed subspace of $\ell^2(\calH)$ for all $N\in\N$.
\begin{proof}
 For a Cauchy sequence $(f^{(n)})_{n\in\N}$ in $\ell^2(\calH)_b^{(N)}$, as above there exist $f_1^{(n)},\dots,f_N^{(n)}\in\calH$ such that we have $f^{(n)}=\sum_{k=1}^Nf_k^{(n)}\bullet b_k$ for $n\in\N$. Hence
\[ \big\|f^{(n)}-f^{(m)}\big\|^2=\left\|\sum_{k=1}^N\big(f_k^{(n)}-f_k^{(m)}\big)\bullet b_k\right\|^2=\sum_{k=1}^N\big\|f_k^{(n)}-f_k^{(m)}\big\|^2\quad\text{for }n,m\in\N, \]
which yields that for all $k=1,\dots,N$ the sequence $\big(f_k^{(n)}\big)_{n\in\N}$ is a Cauchy sequence in $\calH$ with some limit $f_k\in\calH$. Then $f:=\sum_{k=1}^Nf_k\bullet b_k\in\ell^2(\calH)_b^{(N)}$ is the limit of $(f^{(n)})_{n\in\N}$, since
\[ \lim_{n\to\infty}\big\|f-f^{(n)}\big\|^2=\lim_{n\to\infty}\sum_{k=1}^N\big\|f_k-f_k^{(n)}\big\|^2=0. \]
\end{proof}
\end{corollary}

\begin{lemma}\label{someorthorelation}
 Let $f\in\ell^2(\calH)$ and $N\in\N$. Then
$f-\sum_{k=1}^N[f,b_k]\bullet b_k\in\ell^2(\calH)_b^{{(N)}\perp}$.
\begin{proof}
 Let $h_1,\dots,h_N\in\calH$ be arbitrary. By Proposition~\ref{rechenregeln} we have 
\begin{align*}
 \left(f-\sum_{k=1}^N[f,b_k]\bullet b_k,\sum_{k=1}^Nh_k\bullet b_k\right)
 & = \left(f,\sum_{k=1}^Nh_k\bullet b_k\right)-\left(\sum_{k=1}^N[f,b_k]\bullet b_k,\sum_{k=1}^Nh_k\bullet b_k\right)\\
 & = \sum_{k=1}^N(f,h_k\bullet b_k)-\sum_{k=1}^N\big([f,b_k],h_k\big)\\
 & = \sum_{k=1}^N\big([f,b_k],h_k\big)-\sum_{k=1}^N\big([f,b_k],h_k\big)\\
 & = 0.
\end{align*}
\end{proof}
\end{lemma}

\begin{lemma}
 The subspace
\[ \ell^2(\calH)_b^{\infty}:=\bigcup_{N\in\N}\ell^2(\calH)_b^{(N)} \]
is dense in $\ell^2(\calH)$.
\begin{proof}
 Since the subspace of finite sequences is dense in $\ell^2(\calH)$, it suffices to approximate $h\bullet e_s$ for some given $h\in\calH$ and $s\in\N$. We may assume $h\not=0$. Since $(b_k)_{k\in\N}$ is an orthonormal basis of $\ell^2(\R)$ we have $e_s=\sum_{k=1}^\infty(e_s,b_k)b_k$, hence for $\epsilon>0$ there exists $N\in\N$ such that $\big\|e_s-\sum_{k=1}^N(e_s,b_k)b_k\big\|<\epsilon\|h\|^{-1}$. We compute
\[ \left\|h\bullet e_s-\sum_{k=1}^N(e_s,b_k)h\bullet  b_k\right\|=\left\|h\bullet e_s-h\bullet\sum_{k=1}^N(e_s,b_k) b_k\right\|=\|h\|\cdot\left\|e_s-\sum_{k=1}^N(e_s,b_k) b_k\right\|<\epsilon. \]
\end{proof}
\end{lemma}

\begin{remark}[Theorem of best approximation]
 Let $H$ be a pre-Hilbert space and $G$ be a complete subspace of $H$. For a fixed $h\in H$ there exists a unique $g\in G$ such that $\|h-g\|=\dist(h,G):=\inf_{u\in G}\|h-u\|$. Furthermore $g$ is characterized by $h-g\in G^\perp$.
\end{remark}

\begin{theorem}
 For $f\in\ell^2(\calH)$ we have the identity
\[ f=\sum_{k=1}^\infty[f,b_k]\bullet b_k. \]
In particular $\|f\|^2=\sum_{k=1}^\infty\|[f,b_k]\|^2$.
\begin{proof}
 By Lemma~\ref{someorthorelation} for all $f\in\ell^2(\calH)$ and $N\in\N$ it holds
\[ f-\sum_{k=1}^N[f,b_k]\bullet b_k\in\ell^2(\calH)_b^{{(N)}\perp}, \]
and since $\ell^2(\calH)_b^{(N)}$ is closed by Corollary~\ref{finseqclosed}, the theorem of best approximation yields
\[ \left\|f-\sum_{k=1}^N[f,b_k]\bullet b_k\right\|\le\|f-g\|\quad\text{for all }g\in\ell^2(\calH)_b^{(N)}. \]
To a given $\epsilon>0$, we choose $N_0\in\N$ and $g\in\ell^2(\calH)_b^{(N_0)}$ with $\|f-g\|<\epsilon$. Note that $g\in\ell^2(\calH)_b^{(N)}$ for all $N\ge N_0$. We obtain
\[ \left\|f-\sum_{k=1}^N[f,b_k]\bullet b_k\right\|\le\|f-g\|<\epsilon\quad\text{for all }N\ge N_0, \]
i.e. $f=\lim_{N\to\infty}\sum_{k=1}^N[f,b_k]\bullet b_k$. This also yields
\[ \|f\|^2=\lim_{N\to\infty}\left\|\sum_{k=1}^N[f,b_k]\bullet b_k\right\|^2=\lim_{N\to\infty}\sum_{k=1}^N\big\|[f,b_k]\big\|^2. \]
\end{proof}
\end{theorem}

Using this theorem, we can easily generalize Proposition~\ref{finseqdense}:
\begin{corollary}
 If $(h_i)_{i\in\N}$ is an orthonormal basis of $\calH$, then $\{h_i\bullet b_k:i,k\in\N\}$ is an orthonormal basis of $\ell^2(\calH)$.
\end{corollary}

\end{section}

\begin{section}{Operators with Matrix Representation}

\begin{theorem}\label{matrixextension}
 Let $A\in L(\ell^2(\R))$. Then by abuse of notation we define
\begin{equation}\label{extensionformula}
 \ell^2(\calH)\ni f\longmapsto Af:=\sum_{k=1}^\infty[f,b_k]\bullet Ab_k\in\ell^2(\calH)
\end{equation}
and $A$ becomes a bounded linear operator on $\ell^2(\calH)$ with $\|A\|_{L(\ell^2(\calH))}=\|A\|_{L(\ell^2(\R))}$. This definition does not depend on the particular choice of the orthonormal basis $b=(b_k)_{k\in\N}$.
\begin{proof}
We first show that $A$ is a bounded linear operator on the dense subspace $\ell^2(\calH)_b^\infty$ and uniquely extend it to an element of $L(\ell^2(\calH))$. To this end, let $(h_i)_{i\in\N}$ be an orthonormal basis of $\calH$ and $N\in\N$. For $f\in\ell^2(\calH)_b^{(N)}$, $k=1,\dots,N$ and $i\in\N$ we denote
\[ f_{ki}:=\big([f,b_k],h_i\big)_\calH\in\R\quad\text{and}\quad x_i^N:=\sum_{k=1}^Nf_{ki}b_k\in\ell^2(\R). \]
Since $(h_i)_{i\in\N}$ is an orthonormal basis of $\calH$ we have
\begin{equation}\label{fkifliequation}
 \big([f,b_k],[f,b_l]\big)=\sum_{i=1}^\infty\big([f,b_k],h_i\big)\big(h_i,[f,b_l]\big)=\sum_{i=1}^\infty f_{ki}f_{li}
\end{equation}
and thus in particular
\begin{equation}\label{xinequation}
 \sum_{i=1}^\infty\left\|x_i^N\right\|^2=\sum_{i=1}^\infty\sum_{k=1}^Nf_{ki}^2=\sum_{k=1}^N\sum_{i=1}^\infty f_{ki}^2\stackrel{\eqref{fkifliequation}}{=}\sum_{k=1}^N\big\|[f,b_k]\big\|^2=\|f\|^2.
\end{equation}
Using these equations we obtain
\begin{align*}
 \left\|\sum_{k=1}^N[f,b_k]\bullet Ab_k\right\|^2
 & = \sum_{k,l=1}^N\big([f,b_k],[f,b_l]\big)(Ab_k,Ab_l)\\
 & \stackrel{\eqref{fkifliequation}}{=} \sum_{k,l=1}^N\sum_{i=1}^\infty f_{ki}f_{li}(Ab_k,Ab_l)\\
 & = \sum_{i=1}^\infty \left(Ax_i^N,Ax_i^N\right)\\
 & \le \|A\|_{L(\ell^2(\R))}^2\sum_{i=1}^\infty\left\|x_i^N\right\|^2\\
 & \stackrel{\eqref{xinequation}}{=} \|A\|_{L(\ell^2(\R))}^2\|f\|^2,
\end{align*}
which gives rise to $A\in L(\ell^2(\calH))$ with $\|A\|_{L(\ell^2(\calH))}\le\|A\|_{L(\ell^2(\R))}$. For an arbitrary $f\in\ell^2(\calH)$ we then have
\[ Af=A\left(\lim_{N\to\infty}\sum_{k=1}^N[f,b_k]\bullet b_k\right)=\lim_{N\to\infty}A\sum_{k=1}^N[f,b_k]\bullet b_k=\lim_{N\to\infty}\sum_{k=1}^N[f,b_k]\bullet Ab_k, \]
so our definition in \eqref{extensionformula} makes sense. In order to show equality for the operator norms, we note that for $h\in\calH$ and $x\in\ell^2(\calH)$ by Proposition~\ref{rechenregeln} and continuity of $h\bullet\cdot$ it holds
\begin{equation}\label{Ahx=hAx}
 A(h\bullet x)=\sum_{k=1}^\infty[h\bullet x,b_k] \bullet Ab_k=\sum_{k=1}^\infty(x,b_k)h \bullet Ab_k=h\bullet A\left(\sum_{k=1}^\infty(x,b_k)b_k\right)=h\bullet Ax.
\end{equation}
Hence to a given $\epsilon>0$ we choose $x\in\ell^2(\R)$ with $\|Ax\|>(\|A\|_{L(\ell^2(\R))}-\epsilon)\|x\|$ to obtain
\[ \|A(h\bullet x)\|=\|h\bullet Ax\|=\|h\|\cdot\|Ax\|\ge\|h\|(\|A\|_{L(\ell^2(\R))}-\epsilon)\|x\|=(\|A\|_{L(\ell^2(\R))}-\epsilon)\|h\bullet x\|. \]
This gives $\|A\|_{L(\ell^2(\calH))}=\|A\|_{L(\ell^2(\R))}$. Finally, for an arbitrary orthonormal basis $(\beta_k)_{k\in\N}$ of $\ell^2(\R)$ and $f\in\ell^2(\calH)$, by continuity of $A$ and Equation~\eqref{Ahx=hAx} we have
\[ Af=A\left(\sum_{k=1}^\infty[f,\beta_k]\bullet\beta_k\right)=\sum_{k=1}^\infty A\big([f,\beta_k]\bullet\beta_k\big)=\sum_{k=1}^\infty[f,\beta_k]\bullet A\beta_k. \]
\end{proof}
\end{theorem}

\begin{remark}
 If $\lim_{n\to\infty} x_n=x$ and $\lim_{n\to\infty} y_n=y$ in $\calH$, then $\lim_{n\to\infty}(x_n,y_n)=(x,y)$.
\begin{proof}
By the triangle inequality for the modulus and the Cauchy-Bunyakovsky-Schwarz inequality for $(\cdot,\cdot)$ we have
\[ |(x,y)-(x_n,y_n)|\le |(x,y)-(x_n,y)|+|(x_n,y)-(x_n,y_n)|\le \|x-x_n\|\cdot\|y\|+\|x_n\|\cdot\|y-y_n\|, \]
where the right hand side converges to zero since $(x_n)_{n\in\N}$ is bounded.
\end{proof}
\end{remark}

\begin{lemma}
 If $A\in L(\ell^2(\R))$ is (self-adjoint/positive definite), then $A\in L(\ell^2(\calH))$ is (self-adjoint/positive definite).
\begin{proof}
 Let $f,g\in\ell^2(\calH)$. If $A\in L(\ell^2(\R))$ is self-adjoint, then by the remark above we have
\begin{align*}
(f,Ag)_{\ell^2(\calH)}
 & = \lim_{N\to\infty}\left(\sum_{k=1}^N[f,b_k]\bullet b_k,\sum_{l=1}^N[g,b_l]\bullet Ab_l\right)_{\ell^2(\calH)}\\
 & = \lim_{N\to\infty}\sum_{k,l=1}^N\big([f,b_k],[g,b_l]\big)_\calH(b_k,Ab_l)_{\ell^2(\R)}\\
 & = \lim_{N\to\infty}\sum_{k,l=1}^N\big([f,b_k],[g,b_l]\big)_\calH(Ab_k,b_l)_{\ell^2(\R)}\\
 & = \lim_{N\to\infty}\left(\sum_{k=1}^N[f,b_k]\bullet Ab_k,\sum_{l=1}^N[g,b_l]\bullet b_l\right)_{\ell^2(\calH)}\\
 & = (Af,g)_{\ell^2(\calH)}.
\end{align*}
Now assume $A\in L(\ell^2(\R))$ to be positive definite and let $(h_i)_{i\in\N}$ be an orthonormal basis of $\calH$. We abbreviate $f_{ki}:=\big([f,b_k],h_i\big)_\calH\in\R$ for $k,i\in\N$ and set $x_i^N:=\sum_{k=1}^Nf_{ki}b_k\in\ell^2(\R)$ for $N\in\N$ as in the proof of Theorem~\ref{matrixextension}. For $i\in\N$ we then have
\[ \sum_{k=1}^\infty f_{ki}^2=\sum_{k=1}^\infty\big([f,b_k],h_i\big)^2\le\sum_{k=1}^\infty\big\|[f,b_k]\big\|^2\|h_i\|^2=\sum_{k=1}^\infty\big\|[f,b_k]\big\|^2=\|f\|^2<\infty, \]
hence we can define
\[ x_i:=\lim_{N\to\infty}x_i^N=\sum_{k=1}^\infty f_{ki}b_k\in\ell^2(\R) \] 
and obtain $(x_i,Ax_i)_{\ell^2(\R)}\ge 0$ by assumption. Together with Fatou's Lemma we compute
\begin{align*}
(f,Af)_{\ell^2(\calH)}
 & = \lim_{N\to\infty}\left(\sum_{k=1}^N[f,b_k]\bullet b_k,\sum_{l=1}^N[f,b_l]\bullet Ab_l\right)_{\ell^2(\calH)}\\
 & = \lim_{N\to\infty}\sum_{k,l=1}^N\big([f,b_k],[f,b_l]\big)_\calH(b_k,Ab_l)_{\ell^2(\R)}\\
 & \stackrel{\eqref{fkifliequation}}{=} \lim_{N\to\infty}\sum_{k,l=1}^N\sum_{i=1}^\infty f_{ki}f_{li}(b_k,Ab_l)_{\ell^2(\R)}\\
 & = \lim_{N\to\infty}\sum_{i=1}^\infty\big(x_i^N,Ax_i^N\big)_{\ell^2(\R)}\\
 & \ge \sum_{i=1}^\infty\lim_{N\to\infty}\big(x_i^N,Ax_i^N\big)_{\ell^2(\R)}\\
 & = \sum_{i=1}^\infty\underbrace{\big(x_i,Ax_i\big)_{\ell^2(\R)}}_{\ge 0}\\
 & \ge 0.
\end{align*}
If $(f,Af)_{\ell^2(\calH)}=0$, then $(x_i,Ax_i)_{\ell^2(\R)}=0$ for all $i\in\N$ and therefore $f_{ki}=0$ for all $i,k\in\N$, hence $[f,b_k]=0$ for all $k\in\N$ which yields $f=0$.
\end{proof}
\end{lemma}

\begin{corollary}
 If $A\in L(\ell^2(\R))$ is self-adjoint and positive definite, then
\[ \ell^2(\calH)\times\ell^2(\calH)\ni(f,g)\longmapsto(f,g)_A:=(f,Ag)_{\ell^2(\calH)}\in\R \]
defines an inner product on $\ell^2(\calH)$ with corresponding norm $\|\cdot\|_A:=\sqrt{(\cdot,\cdot)_A}$ for which we have $\|\cdot\|_A\le\|A\|^\frac{1}{2}\|\cdot\|$.
\end{corollary}

\begin{corollary}\label{miniregel}
 For $f\in\ell^2(\calH)$, $g,h\in\calH$, $x,y\in\ell^2(\R)$ and a self-adjoint and positive definite operator $A\in L(\ell^2(\R))$ we have the following identities:
\begin{enumerate}
 \item $(h\bullet x,g\bullet y)_A=(h,g)(x,y)_A$.
 \item $(f,h\bullet x)_A=\big([f,Ax],h\big)$.
 \item $[h\bullet x,Ay]=(x,y)_Ah$.
\end{enumerate}
\begin{proof}
 These are direct implications of Proposition~\ref{rechenregeln} together with Equation~\eqref{Ahx=hAx}:
\begin{enumerate}
 \item $(h\bullet x,g\bullet y)_A=(h\bullet x,g\bullet Ay)=(h,g)(x,Ay)=(h,g)(x,y)_A$.
 \item $(f,h\bullet x)_A=(f,h\bullet Ax)=\big([f,Ax],h\big)$.
 \item $[h\bullet x,Ay]=(x,Ay)h=(x,y)_Ah$.
\end{enumerate}
\end{proof}
\end{corollary}

\end{section}

\end{chapter}

\begin{chapter}{The Correlated Gaussian Measure}

\begin{section}{Construction}

Countably Hilbert spaces and in particular nuclear spaces have widely been studied in various literature. We briefly state the following definition for a \emph{Gel'fand triple}, also known as \emph{rigged Hilbert space}, see e.g. \cite{GelVil1964}, which serves our intention to construct a Gaussian measure by means of the Bochner-Minlos theorem. Within the definition we collect some common facts.

\begin{definition}\label{defgelfand}
 Let $\calN$ be a topological vector space and $\calN'$ its topological dual space. We call $\calN\subset\calH\subset\calN'$ a \emph{Gel'fand triple} if the following holds: The topology on $\calN$ is defined by a family of inner products $((\cdot,\cdot)_p)_{p\in\N_0}$ with corresponding norms $(\|\cdot\|_p)_{p\in\N_0}$, which we assume to be compatible in the sense that if $p,q\in\N_0$ and a sequence $(\xi_n)_{n\in\N}$ in $\calN$ converges to zero with respect to $\|\cdot\|_p$ and is a Cauchy sequence with respect to $\|\cdot\|_q$, then $\lim_{n\to\infty}\|\xi_n\|_q=0$. It is easy to show that the topology is induced by the translation invariant metric
\[ \calN\times\calN\ni(\xi,\zeta)\longmapsto\sum_{p=0}^\infty2^{-p}\frac{\|\xi-\zeta\|_p}{1+\|\xi-\zeta\|_p}\in\R. \]
Furthermore assume that $\calN$ is complete with respect to this metric. Without loss of generality we may assume $(\cdot,\cdot)_p\le(\cdot,\cdot)_{p+1}$ for $p\in\N_0$, since otherwise we may replace the family of inner products with the family given by $(\cdot,\cdot)_p':=\sum_{k=0}^p(\cdot,\cdot)_k$, which does not alter the topology on $\calN$ but is monotonously increasing. For $p\in\N_0$ let $\calN_p$ be the Hilbert space obtained by taking the abstract completion of $\calN$ with respect to $\|\cdot\|_p$. Assume $\calN_0=\calH$, which implies $\calN\subset\calH$ densely and continuously. Since the family of norms is increasing, by identifying $\calH$ with its topological dual space $\calH'$, we obtain the chain of spaces
\[ \calN\subset\cdots\subset\calN_2\subset\calN_1\subset\calH=\calH'\subset\calN_{-1}\subset\calN_{-2}\subset\cdots\subset\calN', \]
where $\calN_{-p}$ is the topological dual space of $\calN_p$ for $p\in\N_0$. The completeness of $\calN$ is actually equivalent to $\calN=\bigcap_{p\in\N_0}\calN_p$. It can be shown $\calN'=\bigcup_{p\in\N_0}\calN_{-p}$ and we consider the finest topology on $\calN'$ such that all inclusions $\calN_{-p}\hookrightarrow\calN'$ are continuous. The final important assumption is that for each $p\in\N_0$ the inclusion $N_{p+1,p}:\calN_{p+1}\hookrightarrow\calN_p$ is a Hilbert-Schmidt operator, i.e. for some orthonormal basis $(\eta_k)_{k\in\N}$ in $\calN_{p+1}$ we have
\[ \|N_{p+1,p}\|_{\text{HS}}^2:=\sum_{k=1}^\infty\|\eta_k\|_p^2<\infty, \]
whose value does not depend on the particular choice of the orthonormal basis $(\eta_k)_{k\in\N}$.
\end{definition}

\begin{definition}
 For $p\in\Z$ and some Hilbert space $\big(H,(\cdot,\cdot)_H\big)$ we denote
\[ \ell_p^2(H):=\left\{f\in H^\N:\sum_{k=1}^\infty k^{2p}\|f_k\|_H^2<\infty\right\}, \]
which becomes a Hilbert space itself in an obvious way.
\end{definition}

The following theorem yields a Gef'fand triple with central Hilbert space $\ell^2(\calH)$, if a Gel'fand triple with $\calH$ as central Hilbert space is already given.

\begin{theorem}\label{seqgelfand}
 Assume we have a Gel'fand triple $\calN\subset\calH\subset\calN'$ and let the spaces $\calN_p$, $p\in\N_0$ be as in Definition~\ref{defgelfand}. Then we obtain a Gel'fand triple $s(\calN)\subset\ell^2(\calH)\subset s'(\calN)$ by defining
\[ s(\calN):=\bigcap_{p\in\N_0}\ell_p^2(\calN_p). \]
The topology on $s(\calN)$ we define to be given by the family of norms on $\ell_p^2(\calN_p)$ for $p\in\N_0$.
\begin{proof}
Let $\|\cdot\|_{p,p}$ denote the norm on $\ell_p^2(\calN_p)$ for $p\in\N_0$. Clearly these norms are compatible in the sense of Definition~\ref{defgelfand} and monotonously increasing. Furthermore $\ell_0^2(\calN_0)=\ell^2(\calH)$ and for $p\in\N_0$ the abstract completion of $s(\calN)$ with respect to $\|\cdot\|_{p,p}$ yields exactly $\ell_p^2(\calN_p)$, since $s(\calN)$ contains the set of finite sequences in $\calN$ which are dense in the complete space $\ell_p^2(\calN_p)$. Using $s(\calN)=\bigcap_{p\in\N_0}\ell_p^2(\calN_p)$ yields that the metric
\[ s(\calN)\times s(\calN)\ni(\phi,\psi)\longmapsto\sum_{p=0}^\infty2^{-p}\frac{\|\phi-\psi\|_{p,p}}{1+\|\phi-\psi\|_{p,p}}\in\R \]
is complete, and it clearly induces the topology on $s(\calN)$. It remains to show that for each $p\in\N_0$ the inclusion $I_{p+1,p}:\ell_{p+1}^2(\calN_{p+1})\hookrightarrow\ell_p^2(\calN_p)$ is a Hilbert-Schmidt operator. To this end let $p\in\N_0$ and $(\eta_k)_{k\in\N}$ be an orthonormal basis of $\calN_{p+1}$. Then clearly the set
\[ \left\{\frac{(\eta_k\delta_{l,m})_{m\in\N}}{\|(\eta_k\delta_{l,m})_{m\in\N}\|_{p+1,p+1}}:k,l\in\N\right\}=\left\{\left(\frac{\eta_k\delta_{l,m}}{l^{p+1}}\right)_{m\in\N}:k,l\in\N\right\} \]
is an orthonormal basis of $\ell_{p+1}^2(\calN_{p+1})$. We compute
\[ \|I_{p+1,p}\|_\text{HS}^2=\sum_{k,l=1}^\infty\frac{\|(\eta_k\delta_{l,m})_{m\in\N}\|_{p,p}^2}{\|(\eta_k\delta_{l,m})_{m\in\N}\|_{p+1,p+1}^2}=\sum_{k,l=1}^\infty\|\eta_k\|_p^2\frac{l^{2p}}{l^{2p+2}}=\|N_{p+1,p}\|_\text{HS}^2\sum_{l=1}^\infty\frac{1}{l^2}<\infty, \]
where $N_{p+1,p}:\calN_{p+1}\hookrightarrow\calN_p$ is the inclusion and $\|\cdot\|_p$ denotes the norm on $\calN_p$. This completes the proof.
\end{proof}
\end{theorem}

\begin{example}
Consider the \emph{Schwartz space of functions of rapid decrease}, defined by
\[ S(\R):=\left\{\eta\in C^\infty(\R):\|\eta\|_{n,m}:=\sup_{x\in\R}\left|x^mD^n\eta(x)\right|<\infty\text{ for all }n,m\in\N_0\right\} \]
and equipped with the topology given by the family of seminorms $(\|\cdot\|_{n,m})_{n,m\in\N}$. It is well-known that this is a completely metrizable dense nuclear subspace of the Hilbert space $L^2(\R,\dx)$ and thus yields a Gel'fand triple $S(\R)\subset L^2(\R,\dx)\subset S'(\R)$, which is the standard triple used in White Noise Analysis, see \cite{HKPS1993,ReedSimon1980}. The above theorem can be applied to obtain a Gel'fand triple $s(S(\R))\subset\ell^2(L^2(\R,\dx))\subset s'(S(\R))$.
\end{example}

\begin{notation}
 We denote the canonical dual pairing between $s(\calN)$ and $s'(\calN)$ by
\[ s(\calN)\times s'(\calN)\ni(\phi,\omega)\longmapsto\la\phi,\omega\ra:=\omega(\phi)\in\R. \]
Since we identify $\ell^2(\calH)$ with its dual, for $\phi\in s(\calN)$ and $\omega\in\ell^2(\calH)\subset s'(\calN)$ we have
\[ \la\phi,\omega\ra=(\phi,\omega)_{\ell^2(\calH)}. \]
\end{notation}

\begin{definition}
 We equip $s'(\calN)$ with the $\sigma$-algebra generated by the mappings
\[ s'(\calN)\ni\omega\longmapsto\big(\la\phi_1,\omega\ra,\dots,\la\phi_n,\omega\ra\big)\in\R^n,\quad\text{for }n\in\N\text{ and }\phi_1,\dots,\phi_n\in s(\calN), \]
which is also called the \emph{cylindrical} $\sigma$-algebra.
\end{definition}

The Bochner-Minlos theorem is the standard tool used to obtain a Gaussian measure on spaces like $s'(\calN)$, see \cite{Obata1994}.

\begin{theorem}[Bochner-Minlos theorem]
 Let $C:s(\calN)\to\C$ be a characteristic function, in other words we have $C(0)=1$ and $C$ is continuous and positive semidefinite, i.e.
\[ \sum_{i,j=1}^n\alpha_i\overline{\alpha_j}C(\phi_i-\phi_j)\ge 0\quad\text{ for all }n\in\N\text{ and }\alpha_i\in\C,\,\phi_i\in\calN\text{ for }i=1,\dots,n. \]
Then there exists a unique measure $\mu$ on $s'(\calN)$ which fulfills
\[ \int_{s'(\calN)}\exp\big(i\la\phi,\omega\ra\big)d\mu(\omega)=C(\phi)\quad\text{for all }\phi\in s(\calN). \]
Clearly the measure obtained is a probability measure, since
\[ \mu(s'(\calN))=\int_{s'(\calN)}1d\mu(\omega)=\int_{s'(\calN)}\exp(i\la0,\omega\ra)d\mu(\omega)=C(0)=1. \]
\end{theorem}

\begin{theorem}
  Let $(\cdot,\cdot)'$ be any inner product on $\ell^2(\calH)$ which is continuous. Then
\[ s(\calN)\ni\phi\longmapsto C(\phi):=\exp\left(-\frac{1}{2}(\phi,\phi)'\right)\in\C \]
is a characteristic function in the sense of the Bochner-Minlos theorem.
\begin{proof}
 The equality $C(0)=1$ is clear. Furthermore $C$ is continuous since the embedding $s(\calN)\subset\ell^2(\calH)$ is continuous and $(\cdot,\cdot)'$ is continuous on $\ell^2(\calH)$. Let $n\in\N$ and $\phi_i\in\calN$ for $i=1,\dots,n$. Due to the fact
\[ \sum_{i,j=1}^n\alpha_i\alpha_j(\phi_i,\phi_j)'=\left(\sum_{i=1}^n\alpha_i\phi_i,\sum_{j=1}^n\alpha_j\phi_j\right)'\ge 0\quad\text{for all }\alpha\in\R^n,\]
the matrix $((\phi_i,\phi_j)')_{i,j=1,\dots,n}$ and thus also $(\exp((\phi_i,\phi_j)'))_{i,j=1,\dots,n}$ is positive semidefinite by Lemma~\ref{realposdeflemma} and Corollary~\ref{expisposdef}, see page~\pageref{realposdeflemma}. Now let $\alpha\in\C^n$ be arbitrary. We compute
\begin{align*}
 \sum_{i,j=1}^n\alpha_i\overline{\alpha_j}C(\phi_i-\phi_j)
 & = \sum_{i,j=1}^n\alpha_i\overline{\alpha_j}\exp\left(-\frac{1}{2}\big((\phi_i,\phi_i)'-2(\phi_i,\phi_j)'+(\phi_j,\phi_j)'\big)\right)\\
 & = \sum_{i,j=1}^n\beta_i\overline{\beta_j}\exp\left((\phi_i,\phi_j)'\right)\\
 & \ge 0,
\end{align*}
where $\beta_i=\alpha_i\exp\left(-\frac{1}{2}(\phi_i,\phi_i)'\right)$ for $i=1,\dots,n$. 
\end{proof}
\end{theorem}

\begin{definition}\label{defmuA}
 Let $A\in L(\ell^2(\calH))$ be self-adjoint and positive definite and denote the inner product it generates on $\ell^2(\calH)$ by
\[ \ell^2(\calH)\times\ell^2(\calH)\ni(g,h)\longmapsto(g,h)_A:=(g,Ah)_{\ell^2(\calH)} \]
with corresponding norm $\|\cdot\|_A:=\sqrt{(\cdot,\cdot)_A}$. Since $A$ is continuous, so is $(\cdot,\cdot)_A$. The unique measure $\mu_A$ on $s'(\calN)$ fulfilling 
\[ \int_{s'(\calN)}\exp\left(i\la\phi,\omega\ra\right)d\mu_A(\omega)=\exp\left(-\frac{1}{2}(\phi,\phi)_A\right)\quad\text{for all }\phi\in s(\calN), \]
which exists due to the above theorem, we call the \emph{Gaussian measure with covariance operator $A$}. We denote $L^2(\mu_A):=L^2(s'(\calN),\mu_A;\R)$ and by abuse of notation we will denote the norm on $L^2(\mu_A)$ again by $\|\cdot\|_A$. While Gaussian analysis is usually performed on the complexification of this space, we will stick to the real setting as it suffices for our purposes. However, the results may be transferred to the complex case. To save some space in our equations, we will simply write $s'$ instead of $s'(\calN)$ when integrating, so $\int_{s'(\calN)}fd\mu_A$ becomes $\int_{s'}fd\mu_A$ for integrable or non-negative measureable $f$.
\end{definition}

For the rest of this thesis, $A\in L(\ell^2(\calH))$ will assumed to be self-adjoint and positive definite.

\end{section}

\begin{section}{Properties}

\begin{remark}\label{imagemeasure}
 If $(\Omega,\calF,m)$ is a measure space, $(\Omega',\calF')$ a measureable space and $T:\Omega\to\Omega'$ a measureable map, then $T(m):=m\circ T^{-1}$ is a measure on $\Omega'$, called the \emph{image measure of $m$ under $T$}, and for any measureable $f:\Omega'\to\R$ which is either integrable or non-negative we have
\begin{equation}\label{imagemeasureformula}
 \int_{\Omega'}f(\omega')dT(m)(\omega')=\int_\Omega f(T(\omega))dm(\omega)
\end{equation}
in the sense that either both sides are infinite or both sides are finite and take the same value. Clearly if $T=T'$ almost surely for some for measureable $T':\Omega\to\Omega'$, then $T(m)=T'(m)$.
\end{remark}

\begin{definition}
 By $\mu_n$ we denote the standard Gaussian measure on the measureable space $(\R^n,\mathcal{B}(\R^n))$, i.e. the measure defined by
\[ \mu_n(B)=\left(\frac{1}{\sqrt{2\pi}}\right)^n\int_{B}\exp\left(-\frac{1}{2}|x|^2\right)dx\quad\text{for }B\in\mathcal{B}(\R^n). \]
It is uniquely characterized by its Fourier transform
\[ \R^n\ni p\longmapsto\int_{\R^n}\exp\big(i(p,x)_{\R^n}\big) d\mu_n(x)=\exp\Big(-\frac{1}{2}|p|^2\Big)\in\R. \]
\end{definition}

\begin{lemma}\label{orthogauss}
 Let $n\in\N$ and $\phi_1,\dots,\phi_n\in s(\calN)$ be orthonormal with respect to $(\cdot,\cdot)_A$. Then the image measure of $\mu_A$ under
\[ s'(\calN)\ni\omega\longmapsto T(\omega):=\big(\la\phi_1,\omega\ra,\dots,\la\phi_n,\omega\ra\big)\in\R^n \]
is the standard Gaussian measure $\mu_n$ on $\R^n$.
\begin{proof}
 By Formula~\eqref{imagemeasureformula} from Remark~\ref{imagemeasure} for $p\in\R^n$ we have
\begin{align*}
 \int_{\R^n}\exp\big(i(p,x)_{\R^n}\big)dT(\mu_A)(x)
 & = \int_{s'}\exp\Big(i\big(p,T(\omega)\big)_{\R^n}\Big)d\mu_A(\omega)\\
 & = \int_{s'}\exp\bigg(i\sum_{j=1}^m p_j\la\phi_j,\omega\ra\bigg)d\mu_A(\omega)\\
 & = \int_{s'}\exp\bigg(i\Big\la\sum_{j=1}^m p_j\phi_j,\omega\Big\ra\bigg)d\mu_A(\omega)\\
 & = \exp\bigg(-\frac{1}{2}\Big\|\sum_{j=1}^n p_j\phi_j\Big\|_A^2\bigg)\\
 & = \exp\bigg(-\frac{1}{2}\sum_{j=1}^n p_j^2\bigg)\\
 & = \exp\Big(-\frac{1}{2}|p|^2\Big),
\end{align*}
hence $\mu_n$ and $T(\mu_A)$ have the same Fourier transforms, thus $\mu_n=T(\mu_A)$.
\end{proof}
\end{lemma}

\begin{corollary}\label{orthofubini}
 Let $n\in\N$ and $\phi_1,\dots,\phi_n\in s(\calN)$ be orthonormal with respect to $(\cdot,\cdot)_A$. If for each $i=1,\dots,n$ we have that the measureable function $G_i:\R\to\R$ is non-negative or integrable with respect to the Gaussian measure $\mu_1$, then
\[ \int_{s'}\prod_{i=1}^nG_i(\la\phi_i,\omega\ra)d\mu_A(\omega)=\prod_{i=1}^n\int_{s'}G_i(\la\phi_i,\omega\ra)d\mu_A(\omega). \]
\begin{proof}
 Since $\mu_n$ is the product measure of $n$ one-dimensional measures $\mu_1$, Fubini's Theorem and the lemma above yield
\begin{align*}
 \int_{s'}\prod_{i=1}^nG_i(\la\phi_i,\omega\ra)d\mu_A(\omega)
 & = \int_{\R^n}\prod_{i=1}^nG_i(x_i)d\mu_n(x_1,\dots,x_n)\\
 & = \prod_{i=1}^n\int_{\R}G_i(x_i)d\mu_1(x_i)\\
 & = \prod_{i=1}^n\int_{s'}G_i(\la\phi_i,\omega\ra)d\mu_A(\omega).
\end{align*}
\end{proof}
\end{corollary}

The following yields an isometry from $s(\calN)$ to $L^2(\mu_A)$.

\begin{lemma}\label{isoembedding}
 Let $\phi\in s(\calN)$. Then $\la\phi,\cdot\ra\in L^2(\mu_A)$ with $\|\la\phi,\cdot\ra\|_A=\|\phi\|_A$. 
\begin{proof}
 For $\phi=0$ the statement is clear. Otherwise by Lemma~\ref{orthogauss} we have
\[ \|\la\phi,\cdot\ra\|_A^2=\int_{s'}\la\phi,\omega\ra^2d\mu_A(\omega)=\|\phi\|_A^2\int_{s'}\Big\la\frac{\phi}{\|\phi\|_A},\omega\Big\ra^2d\mu_A(\omega)=\|\phi\|_A^2\int_\R x^2d\mu_1(x)=\|\phi\|_A^2, \]
where we used the well-known fact $\int_\R x^2d\mu_1(x)=1$.
\end{proof}
\end{lemma}

\begin{definition}
 We denote the abstract completion of $\ell^2(\calH)$ with respect to $(\cdot,\cdot)_A$ by $\ell_A^2(\calH)$ and also denote its norm and inner product by $\|\cdot\|_A$ and $(\cdot,\cdot)_A$, respectively.
\end{definition}

\begin{corollary}\label{sinellAcontanddense}
 The inclusion $s(\calN)\subset\ell_A^2(\calH)$ is dense.
\begin{proof}
 To a given $\epsilon>0$ and $f\in\ell_A^2(\calH)$ choose $g\in\ell^2(\calH)$ with $\|f-g\|_A<\epsilon$. For this $g$ there exists $\phi\in s(\calN)$ with $\|g-\phi\|_{\ell^2(\calH)}<\epsilon$. Then
\[ \|f-\phi\|_A\le\|f-g\|_A+\|g-\phi\|_A<\epsilon+\|A\|^\frac{1}{2}\|g-\phi\|_{\ell^2(\calH)}<\left(1+\|A\|^\frac{1}{2}\right)\cdot\epsilon. \]
\end{proof}
\end{corollary}

\begin{lemma}\label{Hilbertkernels}
 Let $f\in\ell_A^2(\calH)$. Since $s(\calN)\subset\ell_A^2(\calH)$ is dense, there exists a sequence $(\phi_k)_{k\in\N}$ in $s(\calN)$ such that $\lim_{k\to\infty}\phi_k=f$ in $\ell_A^2(\calH)$. Then $(\la\phi_k,\cdot\ra)_{k\in\N}$ is a Cauchy sequence in $L^2(\mu_A)$, whose limit is independent of the choice of the approximating sequence $(\phi_k)_{k\in\N}$. Hence $\la f,\cdot\ra:=\lim_{k\to\infty}\la\phi_k,\cdot\ra\in L^2(\mu_A)$ can be defined and for $f\in s(\calN)$ this definition coincides with the equivalence class of the pointwisely defined function $\omega\mapsto\la f,\omega\ra$. Furthermore it holds $\|\la f,\cdot\ra\|_A=\|f\|_A$.
\begin{proof}
 By Lemma~\ref{isoembedding} we have that $(\la\phi_k,\cdot\ra)_{k\in\N}$ is a Cauchy sequence in $L^2(\mu_A)$ and hence converges. If $(\psi_k)_{k\in\N}$ is another sequence in $s(\calN)$ approximating $f$, then for all $k\in\N$ we have
\[ \|\la\phi_k,\cdot\ra-\la\psi_k,\cdot\ra\|_A=\|\phi_k-\psi_k\|_A \le \|\phi_k-f\|_A+\|f-\psi_k\|_A, \]
so $\lim_{k\to\infty}\|\la\phi_k,\cdot\ra-\la\psi_k,\cdot\ra\|_A=0$ and the sequences $(\la\phi_k,\cdot\ra)_{k\in\N}$ and $(\la\psi_k,\cdot\ra)_{k\in\N}$ take the same limit, which we denote by $\la f,\cdot\ra$. By continuity of the norm it holds
\[ \|\la f,\cdot\ra\|_A=\lim_{k\to\infty}\|\la\phi_k,\cdot\ra\|_A=\lim_{k\to\infty}\|\phi_k\|_A=\|f\|_A. \]
\end{proof}
\end{lemma}

\begin{corollary}
 For $f,g\in\ell_A^2(\calH)$ we have $(\la f,\cdot\ra,\la g,\cdot\ra)_A=(f,g)_A$.
\begin{proof}
 By the well-known polarization identity we have
\[ (\la f,\cdot\ra,\la g,\cdot\ra)_A=\frac{1}{4}\left(\|\la f+g,\cdot\ra\|_A^2-\|\la f-g,\cdot\ra\|_A^2\right)=\frac{1}{4}\left(\| f+g\|_A^2-\| f-g\|_A^2\right)=(f,g)_A. \]
\end{proof}
\end{corollary}

\begin{notation}
 Let $(\Omega,\calF,\nu)$ be a measure space and $(f_n)_{n\in\N}$ be a sequence of real-valued measureable functions on $\Omega$ such that the measureable set $N:=\Omega\setminus\{\omega:\lim_{n\to\infty}f_n(\omega)\text{ exists}\}$ has measure zero. Then we define the function $\lim_{n\to\infty}f_n:\Omega\to\R$ by
\[ \left(\lim_{n\to\infty}f_n\right)(\omega):=\lim_{n\to\infty}\indi_{\Omega\setminus N}(\omega)f_n(\omega)=\begin{cases}\lim_{n\to\infty}f_n(\omega)&\omega\in\Omega\setminus N\\0&\omega\in N\end{cases},\quad\omega\in\Omega,\]
which is measureable.
\end{notation}

\begin{remark}
 Let $(\Omega,\calF,\nu)$ be a measure space and let $\lim_{n\to\infty}[f_n]=[f]$ in $L^p(\Omega)$ for some $p\in[1,\infty)$. Then there exists a subsequence $(n_k)_{k\in\N}$ such that $\lim_{k\to\infty}f_{n_k}(\omega)=f(\omega)$ for almost all $\omega\in\Omega$ (or almost surely, if $\nu$ is a probability measure).
\end{remark}

\begin{remark}\label{Hoelderprob}
 Let $(\Omega,\calF,P)$ be a probability space and $1\le p\le q<\infty$. By H\"older's inequality we have $L^q(\Omega,P)\subset L^p(\Omega,P)$ and $\|\cdot\|_{L^p}\le\|\cdot\|_{L^q}$ on $L^q(\Omega,P)$.
\end{remark}

\begin{proposition}\label{generalizedcharfn}
 For $f\in\ell_A^2(\calH)$ we have
\[ \int_{s'}\exp\left(i\la f,\cdot\ra\right)d\mu_A=\exp\left(-\frac{1}{2}(f,f)_A\right). \]
\begin{proof}
 Let $(\phi_k)_{k\in\N}$ be a sequence in $s(\calN)$ with $\lim_{k\to\infty}\phi_k=f$ in $\ell_A^2(\calH)$. Then we have $\lim_{k\to\infty}\la\phi_k,\cdot\ra=\la f,\cdot\ra$ in $L^2(\mu_A)$ by definition. We fix some pointwisely defined representative of $\la f,\cdot\ra$ and also denote it by $\la f,\cdot\ra$. By dropping to a subsequence we may assume that we have $\lim_{k\to\infty}\la\phi_k,\cdot\ra=\la f,\cdot\ra$ almost surely. Since for all $k\in\N$ it holds $|\exp\left(i\la\phi_k,\cdot\ra\right)|=1\in L^2(\mu_A)\subset L^1(\mu_A)$, we apply Lebesgue's theorem of dominated convergence to obtain
\begin{align*}
 \int_{s'}\exp\left(i\la f,\omega\ra\right)d\mu_A(\omega)
 & = \lim_{k\to\infty}\int_{s'}\exp\left(i\la\phi_k,\omega\ra\right)d\mu_A(\omega)\\
 & = \lim_{k\to\infty}\exp\left(-\frac{1}{2}(\phi_k,\phi_k)_A\right)\\
 & = \exp\left(-\frac{1}{2}(f,f)_A\right).
\end{align*}
\end{proof}
\end{proposition}

We now may generalize Lemma~\ref{orthogauss} and Corollary~\ref{orthofubini} in the following ways:

\begin{corollary}\label{orthogaussgeneral}
 Let $n\in\N$ and $f_1,\dots,f_n$ be an orthonormal system in $\ell_A^2(\calH)$. Then the image measure of $\mu_A$ under $T:=\big(\la f_1,\cdot\ra,\dots,\la f_n,\cdot\ra\big)$ is the standard Gaussian measure $\mu_n$ on $\R^n$.
\begin{proof}
 The proof works exactly along the same lines as the proof of Lemma~\ref{orthogauss}, where Remark~\ref{imagemeasure} is used to care about the fact that $T$ is only defined up to almost sure equality and Proposition~\ref{generalizedcharfn} is used to express the Fourier transform in terms of the measure's characterisic function as in the proof of Lemma~\ref{orthogauss}.
\end{proof}
\end{corollary}

\begin{corollary}\label{orthofubinigeneral}
 Let $n\in\N$ and $f_1,\dots,f_n$ be an orthonormal system in $\ell_A^2(\calH)$. If for each $i=1,\dots,n$ we have that the measureable function $G_i:\R\to\R$ is non-negative or integrable with respect to the Gaussian measure $\mu_1$ on $\R$, then
\[ \int_{s'}\prod_{i=1}^nG_i(\la f_i,\cdot\ra)d\mu_A=\prod_{i=1}^n\int_{s'}G_i(\la f_i,\cdot\ra)d\mu_A. \]
\end{corollary}

\end{section}

\begin{section}{The Chaos Decomposition}

\begin{notation}
 For two vector spaces $X$ and $Y$ over the same field $\K$, their \emph{algebraic tensor product}, which again is a vector space over $\K$, is denoted by $X\otimes Y$. Up to an isomorphism in the category of linear $\K$-spaces it is uniquely characterized by the following universal property: There exists a bilinear map $B:X\times Y\to X\otimes Y$ such that for any vector space $V$ over $\K$ and any bilinear map $B':X\times Y\to V$ there exists a unique linear map $L:X\otimes Y\to V$ such that
\[ B'(x,y)=L(B(x,y))\quad\text{for all }x\in X,\,y\in Y. \]
It can be shown that such a space always exists. For $x\in X$, $y\in Y$ we will simply write $x\otimes y$ instead of $B(x,y)$. The algebraic tensor product is associative in the sense that if $Z$ is another vector space over $\K$, then the spaces $(X\otimes Y)\otimes Z$ and $X\otimes(Y\otimes Z)$ are isomorphic in the category of linear $\K$-spaces and we will denote both of the spaces by $X\otimes Y\otimes Z$.
\end{notation}

\begin{lemma}
 For two spaces $X$ and $Y$ over the same field $\K$ we have
\[ X\otimes Y=\spann\{x\otimes y:x\in X,y\in Y\}. \]
\begin{proof}
 Let $S:=\spann\{x\otimes y:x\in X,y\in Y\}$. By Zorn's lemma there exists a subspace $T$ of $X\otimes Y$ such that $X\otimes Y=S\oplus T$. Assume $S\subsetneq X\otimes Y$, so $\dim T\not=0$. Then for $L=\id_{X\otimes Y}$ and $L':X\otimes Y\to X\otimes Y$ defined by $L'(s+t)=s$ for $s\in S$ and $t\in T$, we have $L(x\otimes y)=x\otimes y=L'(x\otimes y)$ for all $x\in X$ and $y\in Y$, but $L\not=L'$. This is a contradiction to the universal property of the algebraic tensor product in view of $V=X\otimes Y$ and $B'=B=\otimes$.
\end{proof}
\end{lemma}

\begin{remark}
 Let $v\in X\otimes Y$. We have proven that there exist $m\in\N$ and $\alpha_k\in\K$, $x_k\in X$ and $y_k\in Y$ for $k=1,\dots,m$ such that $v=\sum_{k=1}^m\alpha_kx_k\otimes y_k$. If we set $x_k':=\alpha_kx_k$ for $k=1,\dots,m$, we obtain the easier representation $v=\sum_{k=1}^mx_k'\otimes y_k$.
\end{remark}

\begin{definition}
 Let $X$ be a real or complex vector space and $n\in\N$. We define the \emph{symmetrization} of $x^{(n)}=x_1\otimes\cdots\otimes x_n\in X\otimes\cdots\otimes X$ to be
\[ \wh{x^{(n)}}:=x_1\wh{\otimes}\cdots\wh{\otimes}x_n:=\frac{1}{n!}\sum_{\sigma\in S_n}x_{\sigma(1)}\otimes\cdots\otimes x_{\sigma(n)}. \]
Here $S_n$ stands for the group of permutations on $\{1,\dots,n\}$. One can show that this induces a linear operator on $X\otimes\cdots\otimes X$. For general $x^{(n)}\in X\otimes\cdots\otimes X$, the symmetrization of $\wh{x^{(n)}}$ again yields $\wh{x^{(n)}}$. If $x^{(n)}=\wh{x^{(n)}}$, we call $x^{(n)}$ symmetric and denote the subspace of symmetric elements by
\[ X\wh{\otimes}\cdots\wh{\otimes}X:=\{x^{(n)}\in X\otimes\cdots\otimes X:x^{(n)}=\wh{x^{(n)}}\}. \]
\end{definition}

\begin{notation}
 Let $X$ be a vector space over $\K\in\{\R,\C\}$, $x\in X$ and $n\in\N$. We denote the $n^\text{th}$ tensor power of $x$ by $x^\otn:=x\otimes\cdots\otimes x\in X\otimes\cdots\otimes X$ and set $x^{\otimes 0}=1\in\K$. If $L_1,\dots,L_n:X\to X$ are linear operators, we uniquely define the linear operator $L_1\otimes\cdots\otimes L_n$ on $X\otimes\cdots\otimes X$ by
\[ L_1\otimes\cdots\otimes L_n(x_1\otimes\cdots\otimes x_n):=L_1x_1\otimes\cdots\otimes L_nx_n\quad\text{for }x_1,\dots,x_n\in X \]
and its symmetrization $L_1\wh{\otimes}\cdots\wh{\otimes}L_n$ by
\[ L_1\wh{\otimes}\cdots\wh{\otimes}L_n(x_1\otimes\cdots\otimes x_n):=L_1x_1\wh{\otimes}\cdots\wh{\otimes} L_nx_n\quad\text{for }x_1,\dots,x_n\in X. \]
For a linear operator $L:X\to X$ we introduce the notations 
\[ L^\otn:=L\otimes\cdots\otimes L,\quad L^\wotn:=L\wh{\otimes}\cdots\wh{\otimes}L\quad\text{and}\quad L^{\otimes 0}:=\id_\K. \]
\end{notation}

\begin{corollary}
 For a real or complex vector space $X$ and $n\in\N$ it holds
\[ X\wh{\otimes}\cdots\wh{\otimes}X=\spann\left\{x_1\wh{\otimes}\cdots\wh{\otimes}x_n:x_1,\cdots,x_n\in X\right\}. \]
\begin{proof}
 Let $x^{(n)}\in X\wh{\otimes}\cdots\wh{\otimes}X$. There exist $m\in\N$ and $x_1^k,\dots,x_n^k\in X$ for $k=1,\dots,m$ such that $x^{(n)}=\sum_{k=1}^m x_1^k\otimes\cdots\otimes x_n^k$. Thus
\[ x^{(n)}=\wh{x^{(n)}}=\sum_{k=1}^m x_1^k\wh{\otimes}\cdots\wh{\otimes} x_n^k\in\spann\left\{x_1\wh{\otimes}\cdots\wh{\otimes}x_n:x_1,\cdots,x_n\in X\right\}. \]
\end{proof}
\end{corollary}

The following can be found in \cite{Obata1994}:

\begin{lemma}[Polarization formula]\label{polarizationformula}
 Let $X$ and $Y$ be real or complex vector spaces, $n\in\N$ and $F:X^n\to Y$ be multilinear and symmetric. Then for $x_1,\dots,x_n\in X$ it holds
\[ F(x_1,\dots,x_n)=\frac{1}{2^nn!}\sum_{B\in\{\pm 1\}^n}B_1\cdots B_nA(B_1x_1+\dots+B_nx_n), \]
where $A(x):=F(x,\dots,x)$ for $x\in X$.
\end{lemma}

\begin{corollary}
 For a real or complex space $X$ and $n\in\N$ we have
\[ X\wh{\otimes}\cdots\wh{\otimes}X=\spann\left\{x^\otn:x\in X\right\}. \]
\begin{proof}
 Clearly $X\wh{\otimes}\cdots\wh{\otimes}X$ contains all elements of the form $x^\otn$, where $x\in X$. For the other inclusion define $F:X^n\to X\otimes\cdots\otimes X$ by $F(x_1,\dots,x_n):=x_1\wh{\otimes}\cdots\wh{\otimes} x_n$ for $x_1,\dots,x_n\in X$. Applying the polarization formula yields $x_1\wh{\otimes}\cdots\wh{\otimes} x_n\in\spann\{x^\otn:x\in X\}$. Thus
\[ X\wh{\otimes}\cdots\wh{\otimes}X=\spann\left\{x_1\wh{\otimes}\cdots\wh{\otimes} x_n:x_1,\dots,x_n\in X\right\}=\spann\left\{x^\otn:x\in X\right\}. \]
\end{proof}
\end{corollary}

\begin{definition}
 For $n\in\N$ we denote $s(\calN)^\otn:=s(\calN)\otimes\cdots\otimes s(\calN)$ and its subspace of symmetric elements by $s(\calN)^\wotn:=s(\calN)\wh{\otimes}\cdots\wh{\otimes}s(\calN)$. We set $s(\calN)^{\otimes 0}:=s(\calN)^{\wot 0}:=\R$.
\end{definition}

\begin{definition}
 Let $n\in\N$. By $\ell_A^2(\calH)^\otn$ we denote the abstract completion of the space $\ell_A^2(\calH)\otimes\cdots\otimes\ell_A^2(\calH)$ with respect to the unique inner product which fulfills
\[ \big(f_1\otimes\cdots\otimes f_n,g_1\otimes\cdots\otimes g_n\big)_A:=\prod_{k=1}^n(f_k,g_k)_A\quad\text{for }f_1,\dots,f_n,g_1,\dots,g_n\in\ell_A^2(\calH). \]
The space $\ell_A^2(\calH)^\wotn$ is defined to be the closure of $\ell_A^2(\calH)\wh{\otimes}\cdots\wh{\otimes}\ell_A^2(\calH)$ in $\ell_A^2(\calH)^\otn$, i.e.
\[ \ell_A^2(\calH)^\wotn=\overline{\spann\left\{f^\otn:f\in\ell_A^2(\calH)\right\}}\subset\ell_A^2(\calH)^\otn. \]
Further we set $\ell_A^2(\calH)^{\otimes 0}:=\ell_A^2(\calH)^{\wot 0}:=\R$.
\end{definition}

\begin{remark}
 It is noted that we defined $s(\calN)^\otn$ simply as an algebraic tensor product, while $\ell_A^2(\calH)^\otn$ is the abstract completion of the algebraic tensor product with respect to some inner product. This may seem confusing, but for our purposes this yields the simplest notation.
\end{remark}

\begin{remark}
 If $L_1,\dots,L_n\in L(\ell_A^2(\calH))$, then $L_1\otimes\cdots\otimes L_n\in L\left(\ell_A^2(\calH)\otimes\cdots\otimes\ell_A^2(\calH)\right)$ with operator norm $\|L_1\otimes\cdots\otimes L_n\|=\|L_1\|\cdots\|L_n\|$ and hence can be extended to an element of $L\left(\ell_A^2(\calH)^\otn\right)$, see e.g. \cite{Dixmier1981}.
\end{remark}

\begin{corollary}\label{densetensors} 
 For all $n\in\N$ the inclusions $s(\calN)^\otn\subset\ell_A^2(\calH)^\otn$ and $s(\calN)^\wotn\subset\ell_A^2(\calH)^\wotn$ are dense.
\begin{proof}
 For $n=1$ this is Corollary~\ref{sinellAcontanddense}. Assume the density of $s(\calN)^\otn\subset\ell_A^2(\calH)^\otn$ has been proven for some $n\in\N$ and let $F\in\ell_A^2(\calH)^\otn$ and $f\in\ell_A^2(\calH)$. It suffices to approximate $F\otimes f$, since the linear span of such elements is dense in $\ell_A^2(\calH)^{\otimes n+1}$. We may assume $F\not=0$ and $f\not=0$ since otherwise $F\otimes f=0\in s(\calN)^{\otimes n+1}$. For $\epsilon>0$ choose $\Phi\in s(\calN)^\otn$ with $\|F-\Phi\|_A<\epsilon\|f\|_A^{-1}$. We may enforce $\Phi\not=0$ since $F\not=0$. Choose $\phi\in s(\calN)$ with $\|f-\phi\|_A<\epsilon\|\Phi\|_A^{-1}$. Then
\begin{align*}
 \|F\otimes f-\Phi\otimes\phi\|_A
 & \le \|F\otimes f-\Phi\otimes f\|_A+\|\Phi\otimes f-\Phi\otimes\phi\|_A\\
 & = \|F-\Phi\|_A\|f\|_A+\|\Phi\|_A\|f-\phi\|_A\\
 & < 2\epsilon.
\end{align*}
To prove the density of $s(\calN)^\wotn\subset\ell_A^2(\calH)^\wotn$, we observe that for $f\in\ell_A^2(\calH)$ and a sequence $(\phi_k)_{k\in\N}$ in $s(\calN)$ with $\lim_{k\to\infty}\phi_k=f$ it holds
\begin{align*}
 \lim_{k\to\infty}\|f^\otn-\phi_k^\otn\|_A^2
 & = \lim_{k\to\infty}\|f^\otn\|_A^2+\|\phi_k^\otn\|_A^2-2(f^\otn,\phi_k^\otn)_A\\
 & = \lim_{k\to\infty}\|f\|_A^{2n}+\|\phi_k\|_A^{2n}-2(f,\phi_k)_A^n\\
 & =0.
\end{align*}
Since $\spann\{f^\otn:f\in\ell_A^2(\calH)\}$ is dense in $\ell_A^2(\calH)^\wotn$ by definition, the statement is proven.
\end{proof}
\end{corollary}

\begin{definition}
  The set of monomials $\calM_n$ of order $n\in\N_0$ on $s'(\calN)$ we define by
\[ \calM_n:=\left\{s'(\calN)\ni\omega\mapsto\la\phi^{(n)},\omega^\otn\ra\in\R:\phi^{(n)}\in s(\calN)^\otn\right\}. \]
Since for $\omega\in s'(\calN)$, $n\in\N_0$ and $\phi^{(n)}\in s(\calN)^\otn$ it holds $\la\phi^{(n)},\omega^\otn\ra=\la\wh{\phi^{(n)}},\omega^\otn\ra$, the polarization formula yields
\[ \calM_n=\spann\left\{s'(\calN)\ni\omega\mapsto\la\phi^\otn,\omega^\otn\ra=\la\phi,\omega\ra^n\in\R:\phi\in s(\calN)\right\}. \]
Furthermore, we define the set $\calP_n$ of polynomials of degree $n\in\N_0$ on $s'(\calN)$ and the set of polynomials $\calP$ on $s'(\calN)$ by
\[ \calP_n:=\sum_{k=0}^n\calM_k \quad\text{and}\quad \calP:=\bigcup_{n\in\N_0}\calP_n, \]
respectively. 
\end{definition}

A proof of the following important result works along the same lines as the corresponding proof in \cite{Obata1994}:

\begin{theorem}
 The set of polynomials is dense in $L^2(\mu_A)$.
\end{theorem}

\begin{definition}
 We define $\tau_A:s(\calN)^{\otimes 2}\to\R$ as the unique linear extension of the operator fulfilling $\tau_A(\phi\otimes\psi):=\la\phi,A\psi\ra=(\phi,\psi)_A$ for $\phi,\psi\in s(\calN)$, which exists due to the universal property of the tensor product.
\end{definition}

\begin{definition}
 For $\omega\in s'(\calN)$ and $n\in\N_0$ we inductively define $\colon\omega^\otn\colon\in\big(s(\calN)^\otn\big)^*$, the algebraic dual space of $s(\calN)^\otn$, by
\[ \colon\omega^{\otimes 0}\colon:=\id_\R,\quad\colon\omega^{\otimes 1}\colon:=\omega,\quad\text{and} \]
\[ \colon\omega^\otn\colon:=\omega\wot\colon\omega^{\otimes n-1}\colon-(n-1)\tau_A\wot\colon\omega^{\otimes n-2}\colon\quad\text{for }n\ge2. \]
It is clear from the definition that for $\phi^{(n)}\in s(\calN)^\otn$ we have
\[ \wickpolyo{\phi^{(n)}}=\wickpolyo{\wh{\phi^{(n)}}}. \]
\end{definition}

\begin{lemma}\label{wickpolynomial}
 For $\omega\in s'(\calN)$ and $n\in\N_0$ one has
\[ \colon\omega^\otn\colon=\sum_{k=0}^{\lfloor\frac{n}{2}\rfloor}\frac{(-1)^kn!}{2^kk!(n-2k)!}\tau_A^{\wot k}\wot\omega^{\otimes n-2k}\quad\text{and} \]
\[ \omega^\otn=\sum_{k=0}^{\lfloor\frac{n}{2}\rfloor}\frac{n!}{2^kk!(n-2k)!}\tau_A^{\wot k}\wot\colon\omega^{\otimes n-2k}\colon. \]
\begin{proof}
This proof uses a straightforward but unattractive induction, which will be given in Section~\ref{wickpolynomialproofsection} of the appendix on page~\pageref{wickpolynomialproof}.
\end{proof}
\end{lemma}

\begin{corollary}\label{wickhermite}
 For $\omega\in s'(\calN)$, $\phi\in s(\calN)$, $\phi\not=0$ and $n\in\N_0$ we have
\begin{equation}
 \wickpolyo{\phi^\otn}=\|\phi\|_A^nH_n\left(\frac{\la\phi,\omega\ra}{\|\phi\|_A}\right),
\end{equation}
where $H_n$ is the $n^\text{th}$ Hermite polynomial, see section~\ref{Hermitepolynomial} in the appendix on page~\pageref{Hermitepolynomial}.
\begin{proof}
 This is a direct implication of the previous lemma using Equation~\eqref{Hermitesumrepresentation}.
\end{proof}
\end{corollary}

\begin{definition}
 The set of wick ordered polynomials on $s'(\calN)$ we define by
\[ \calW:=\left\{s'(\calN)\ni\omega\mapsto\sum_{n=0}^m\wickpolyo{\phi^{(n)}}\in\R:m\in\N_0, \phi^{(n)}\in s(\calN)^\otn\text{ for }n=0,\dots,m\right\}. \]
Lemma~\ref{wickpolynomial} yields $\calW=\calP$, hence the wick ordered polynomials are dense in $L^2(\mu_A)$.
\end{definition}

\begin{lemma}
 Let $\phi,\psi\in s(\calN)$ and $n,m\in\N_0$. Then we have
\[ \int_{s'}\wickpolyo{\phi^\otn}\wickpolyo[m]{\psi^\otm}d\mu_A(\omega)=\delta_{n,m}n!(\phi,\psi)_A^n. \]
\begin{proof}
 Without loss of generality we may assume $\|\phi\|_A=\|\psi\|_A=1$. By $\dim\spann\{\phi,\psi\}\le 2$ there exists $\eta\in s(\calN)$ with $\|\eta\|_A=1$ and $(\eta,\psi)_A=0$ such that $\phi\in\spann\{\eta,\psi\}$. Then for $\alpha:=(\phi,\psi)_A$ and $\beta:=(\phi,\eta)_A$ we have $\phi=\alpha\psi+\beta\eta$ with $\alpha^2+\beta^2=\|\phi\|_A=1$. Now for $\omega\in s'(\calN)$, if necessary using the convention $0^0:=1$, Corollary~\ref{wickhermite} above and Equation~\eqref{Hermitebinomial} from page~\pageref{Hermitebinomial} yield
\begin{align}
 \wickpolyo{\phi^\otn}\wickpolyo[m]{\psi^\otm}
 & = H_n(\la\phi,\omega\ra)H_m(\la\psi,\omega\ra)\notag\\
 & = H_n\big(\alpha\la\psi,\omega\ra+\beta\la\eta,\omega\ra\big)H_m(\la\psi,\omega\ra)\label{somelabel1}\\
 & = \sum_{k=0}^n\binom{n}{k}\alpha^k\beta^{n-k}H_k(\la\psi,\omega\ra)H_{n-k}(\la\eta,\omega\ra)H_m(\la\psi,\omega\ra).\notag
\end{align}
By using $\|\eta\|_A=1$, Corollary~\ref{orthofubini} and Equations~\eqref{Hermiteinnerprod} and~\eqref{Hermiteintegraliszero}, for fixed $k\in\{0,\dots,n\}$ we have
\begin{align}
 \int_{s'}H_k(\la\psi,\omega\ra)H_{n-k}(\la\eta,\omega\ra)H_m(\la\psi,\omega\ra)d\mu_A(\omega)
 & = \int_\R H_k(x)H_m(x)d\mu_1(x)\int_\R H_{n-k}(y)d\mu_1(y)\notag\\
 & = (H_k,H_m)_{L^2(\mu_1)}\delta_{n-k,0}\notag\\
 & = m!\delta_{k,m}\delta_{k,n}=n!\delta_{k,n,m}.\label{somelabel2}
\end{align}
Using these equations we obtain
\begin{align*}
 \delta_{n,m}n!(\phi,\psi)_A^n
 & = \delta_{n,m}n!\alpha^n = \sum_{k=0}^n\alpha^kn!\delta_{k,n,m} = \sum_{k=0}^n\binom{n}{k}\alpha^k\beta^{n-k}n!\delta_{k,n,m}\\
 & \stackrel{\eqref{somelabel2}}{=} \sum_{k=0}^n\binom{n}{k}\alpha^k\beta^{n-k}\int_{s'}H_k(\la\psi,\omega\ra)H_{n-k}(\la\eta,\omega\ra)H_m(\la\psi,\omega\ra)d\mu_A(\omega)\\
 & \stackrel{\eqref{somelabel1}}{=} \int_{s'}\wickpolyo{\phi^\otn}\wickpolyo[m]{\psi^\otm}d\mu_A(\omega).
\end{align*}
\end{proof}
\end{lemma}

\begin{corollary}
 Let $n,m\in\N_0$ and $\phi^{(n)}\in s(\calN)^\otn$, $\psi^{(m)}\in s(\calN)^\otm$. Then
\[ \int_{s'}\wickpolyo{\phi^{(n)}}\wickpolyo[m]{\psi^{(m)}}d\mu_A(\omega)=\delta_{n,m}n!\left(\wh{\phi^{(n)}},\wh{\psi^{(m)}}\right)_A. \]
\begin{proof}
 By the polarization formula there exist $r_1,r_2\in\N$, $\phi_1,\dots,\phi_{r_1},\psi_1,\dots\psi_{r_2}\in s(\calN)$ and $\alpha_1,\dots\alpha_{r_1},\beta_1,\dots,\beta_{r_2}\in\R$ such that 
\[ \wh{\phi^{(n)}}=\sum_{k=1}^{r_1}\alpha_i\phi_i^\otn\quad\text{and}\quad\wh{\psi^{(m)}}=\sum_{j=1}^{r_2}\beta_j\psi_j^\otm. \]
Then the previous lemma yields
\begin{align*}
 \int_{s'}\wickpolyo{\phi^{(n)}}\wickpolyo[m]{\psi^{(m)}}d\mu_A(\omega)
 & = \int_{s'}\wickpolyo{\wh{\phi^{(n)}}}\wickpolyo[m]{\wh{\psi^{(m)}}}d\mu_A(\omega)\\
 & = \sum_{i=1}^{r_1}\sum_{j=1}^{r_2}\alpha_i\beta_j\delta_{n,m}n!(\phi_i,\psi_j)_A^n\\
 & = \delta_{n,m}n!\left(\sum_{i=1}^{r_1}\alpha_i\phi_i^\otn,\sum_{j=1}^{r_2}\beta_j\psi_j^\otm\right)_A\\
 & = \delta_{n,m}n!\left(\wh{\phi^{(n)}},\wh{\psi^{(m)}}\right)_A.
\end{align*}
\end{proof}
\end{corollary}

\begin{proposition}
 Let $n\in\N_0$ and $f^{(n)}\in\ell_A^2(\calH)^\wotn$. Similar as in Lemma~\ref{Hilbertkernels} we may define $\wickpoly{f^{(n)}}\in L^2(\mu_A)$ as the element $\lim_{k\to\infty}\wickpoly{\phi_k^{(n)}}$ in $L^2(\mu_A)$, where $(\phi_k^{(n)})_{k\in\N}$ is an arbitrary sequence in $s(\calN)^\wotn$ with $\lim_{k\to\infty}\phi_k^{(n)}=f^{(n)}$ in $\ell_A^2(\calH)^\otn$, whose particular choice is irrelevant. Furthermore we have $\left\|\wickpoly{f^{(n)}}\right\|_A^2=n!\|f^{(n)}\|_A^2$.
\begin{proof}
 Let $f^{(n)}\in\ell_A^2(\calH)^\wotn$. Due to Corollary~\ref{densetensors} there exists a sequence $(\phi_k^{(n)})_{k\in\N}$ in $s(\calN)^\wotn$ with $\lim_{k\to\infty}\phi_k^{(n)}=f^{(n)}$. By the above Corollary we have
\[ \Big\|\wickpoly{\phi_k^{(n)}}\Big\|_A^2=n!\Big\|\wh{\phi_k^{(n)}}\Big\|_A^2=n!\Big\|\phi_k^{(n)}\Big\|_A^2, \]
hence $\left(\wickpoly{\phi_k^{(n)}}\right)_{k\in\N}$ is a Cauchy sequence in $L^2(\mu_A)$. If $(\psi_k^{(n)})_{k\in\N}$ is another sequence in $s(\calN)^\wotn$ with limit $f^{(n)}$, then
\[ \Big\|\wickpoly{\phi_k^{(n)}}-\wickpoly{\psi_k^{(n)}}\Big\|_A^2=n!\Big\|\wh{\phi_k^{(n)}}-\wh{\psi_k^{(n)}}\Big\|_A^2=n!\Big\|\phi_k^{(n)}-\psi_k^{(n)}\Big\|_A^2, \]
so $\lim_{k\to\infty}\left\|\wickpoly{\phi_k^{(n)}}-\wickpoly{\psi_k^{(n)}}\right\|_A=0$ and the sequences $\left(\wickpoly{\phi_k^{(n)}}\right)_{k\in\N}$ and $\left(\wickpoly{\psi_k^{(n)}}\right)_{k\in\N}$ take the same limit, which we denote by $\wickpoly{f^{(n)}}$. By continuity of the norm it holds
\[ \left\|\wickpoly{f^{(n)}}\right\|_A^2=\lim_{k\to\infty}\left\|\wickpoly{\phi_k^{(n)}}\right\|_A^2=\lim_{k\to\infty}n!\left\|\phi_k^{(n)}\right\|_A^2=n!\left\|f^{(n)}\right\|_A^2. \]
\end{proof}
\end{proposition}

This directly implies the following:

\begin{corollary}
 Let $n,m\in\N_0$ and $f^{(n)}\in\ell_A^2(\calH)^\wotn$, $g^{(m)}\in\ell_A^2(\calH)^\wotm$. Then
\[ \int_{s'}\wickpoly{f^{(n)}}\wickpoly[m]{g^{(m)}}d\mu_A=\delta_{n,m}n!(f^{(n)},g^{(m)})_A. \]
\end{corollary}

\begin{theorem}[Chaos decomposition]
 Let $F\in L^2(\mu_A)$. Then for each $n\in\N_0$ there exists a unique $f^{(n)}\in\ell_A^2(\calH)^\wotn$ such that $F=\sum_{n=0}^\infty\wickpoly{f^{(n)}}$ in the $L^2(\mu_A)$ sense. We then have $\|F\|_A^2=\sum_{n=0}^\infty n!\|f^{(n)}\|_A^2$.
\begin{proof}
 Let $(F_k)_{k\in\N}$ be a sequence of wick ordered polynomials with $\lim_{k\to\infty}F_k=F$ in $L^2(\mu_A)$. For each $k\in\N$ there exists $m_k\in\N_0$, $\phi_k^{(n)}\in s(\calN)^\wotn$ for $n=0,\dots,m_k$ such that
\[ F_k=\sum_{n=0}^\infty\wickpoly{\phi_k^{(n)}}, \]
where we set $\phi_k^{(n)}=0$ for $n>m_k$. For fixed $n\in\N_0$, we have
\[ \big\|\phi_k^{(n)}-\phi_l^{(n)}\big\|_A^2\le\sum_{m=0}^\infty m!\big\|\phi_k^{(m)}-\phi_l^{(m)}\big\|_A^2=\|F_k-F_l\|_A^2 \]
for any choices of $k,l\in\N$, hence the sequence $\big(\phi_k^{(n)}\big)_{k\in\N}$ is a Cauchy sequence with some limit $f^{(n)}\in\ell_A^2(\calH)^\wotn$. We define
\[ \tilde{F}_j:=\sum_{n=0}^j\wickpoly{f^{(n)}} \]
for $j\in\N$. For $\epsilon>0$ there exists $l\in\N$ with $\|F-F_l\|_A^2<\epsilon/2$. By having in mind that $\phi_l^{(n)}=0$ for $n>m_l$, we see that for all $i,j\ge m_l$ it holds
\begin{align*}
 \|\tilde{F}_j-\tilde{F}_i\|_A^2
 & = \sum_{n=i+1}^j n!\|f^{(n)}\|_A^2\\
 & = \lim_{k\to\infty}\sum_{n=i+1}^j n!\big\|\phi_k^{(n)}\big\|_A^2\\
 & \le \lim_{k\to\infty}\sum_{n=i+1}^j 2n!\big\|\phi_k^{(n)}-\phi_l^{(n)}\big\|_A^2+\sum_{n=i+1}^j 2n!\big\|\phi_l^{(n)}\big\|_A^2\\
 & \le \lim_{k\to\infty}2\|F_k-F_l\|_A^2\\
 & = 2\|F-F_l\|_A^2\\
 & < \epsilon,
\end{align*}
where we used the estimation $(a+b)^2\le2a^2+2b^2$ for $a,b\in\R$. Hence $(\tilde{F}_j)_{j\in\N}$ is a Cauchy sequence in $L^2(\mu_A)$ and thus we can define
\[ \tilde{F}:=\lim_{j\to\infty}\tilde{F}_j=\sum_{n=0}^\infty\wickpoly{f^{(n)}}\in L^2(\mu_A). \]
It remains to show $F=\tilde{F}$. To this end let $p\in\calP$ be arbitrary with representation
\[ p=\sum_{n=0}^m\wickpoly{\psi^{(n)}}\quad\text{for some }m\in\N_0\text{ and }\psi^{(n)}\in s(\calN)^\otn\text{ for }n=0,\dots,m .\]
Then
\[ (F,p)_A=\lim_{k\to\infty}(F_k,p)_A=\lim_{k\to\infty}\sum_{n=0}^mn!\big(\phi_k^{(n)},\wh{\psi^{(n)}}\big)_A=\sum_{n=0}^mn!\big(f^{(n)},\wh{\psi^{(n)}}\big)_A=(\tilde{F},p)_A \]
and hence $F-\tilde{F}\in\calP^\perp=\{0\}$, i.e. $F=\tilde{F}$. It clearly follows
\[ \|F\|_A^2=\lim_{j\to\infty}\|\tilde{F}_j\|_A^2=\lim_{j\to\infty}\sum_{n=0}^j n!\|f^{(n)}\|_A^2. \]
\end{proof}
\end{theorem}

\end{section}

\end{chapter}

\begin{chapter}{Conditional Expectations}

\begin{section}{Representation}

\begin{definition}
 Let $(\Omega,\calF,P)$ be a probability space and $\calG$ be a sub-$\sigma$-algebra of $\calF$. Let $X:\Omega\to\R$ be a non-negative or integrable random variable. A random variable $Y:\Omega\to\R$ is called \emph{conditional expectation of $X$ given $\calG$}, if $Y$ is $\calG$-measureable and $\E[\indi_GX]=\E[\indi_GY]$ holds for all $G\in\calG$. We denote the set of all conditional expectations of $X$ given $\calG$ by $\E[X|\calG]$. If $Z$ is another random variable we denote $\E[X|Z]:=\E[X|\sigma(Z)]$.
\end{definition}

\begin{remark}\label{augmentationremark}
 Let $(\Omega,\calF,P)$ be a probability space, $\calG$ be a sub-$\sigma$-algebra of $\calF$ and $p\ge 1$. One can show that the conditional expectation defines a contractive operator from $L^p(\Omega,\calF,P)$ onto $L^p(\Omega,\calG,P)$. In particular, the conditional expectation is independent of representatives. Via the isometry
\[ L^p(\Omega,\calG,P)\ni[g]\longmapsto[g]_\calF:=\big\{f:f\text{ is }\calF\text{-measureable and } P(f=g)=1\big\}\in L^p(\Omega,\calF,P) \]
we may consider $L^p(\Omega,\calG,P)$ as a closed subspace of $L^p(\Omega,\calF,P)$. For $[X]\in L^p(\Omega,\calF,P)$ we  especially consider $\E[X|\calG]$ as an element of $L^p(\Omega,\calF,P)$.
\end{remark}

\begin{remark}
 Let $(\Omega,\calF,P)$ be a probability space and $\calG$ be a sub-$\sigma$-algebra of $\calF$. Then for $[X]\in L^2(\Omega,\calF,P)$ we have $E[X|\calG]=\calP_\calG([X])$, where
\[ \calP_\calG:L^2(\Omega,\calF,P)\to L^2(\Omega,\calG,P) \]
is the orthogonal projection. In particular, since the orthogonal projection is continuous, for a Cauchy sequence $([X_n])_{n\in\N}$ in $L^2(\Omega)$ we have $\lim_{n\to\infty}\E\big[[X_n]|\calG\big]=\E\big[\lim_{n\to\infty}[X_n]\big|\calG\big]$ in $L^2(\Omega)$. 
\end{remark}

The well-known factorisation lemma will be very useful in our proofs later on:

\begin{lemma}[Factorisation lemma]\label{factorisationlemma}
 Let $(\Omega,\calF)$ be a measureable space and $Y:\Omega\to\R$ be measureable. If $X:\Omega\to\R$ is $\calG:=\sigma(Y)$-measureable, then there exists a measureable $g:\R\to\R$ such that $X=g(Y)$.
\begin{proof}
 If $X$ is an elementary function, there exists $n\in\N$ and $G_1,\dots,G_n\in\calG$, $\alpha_1,\dots,\alpha_n\in\R$ with $X=\sum_{i=1}^n\alpha_i\indi_{G_i}$. Since $\calG=\sigma(Y)$, there exist Borel sets $B_1,\dots,B_n$ with $G_i=Y^{-1}(B_i)$ for $i=1,\dots,n$. Thus
\[ X=\sum_{i=1}^n\alpha_i\indi_{G_i}=\sum_{i=1}^n\alpha_i\indi_{Y^{-1}(B_i)}=\sum_{i=1}^n\alpha_i\indi_{B_i}(Y)=g(Y)\quad\text{for } g:=\sum_{i=1}^n\alpha_i\indi_{B_i}. \]
If $X$ is non-negative, there exist elementary functions $(X_n)_{n\in\N}$ with $\sup_{n\in\N}X_n(\omega)=X(\omega)$ for all $\omega\in\Omega$. For each $n\in\N$ there exists a measureable $g_n$ with $X_n=g_n(Y)$ as above. Then the pointwisely defined function $g:=\sup_{n\in\N}g_n$ again is a $\calG$-measureable function and we have $X=\sup X_n=\sup g_n(Y)=g(Y)$. For some arbitrary measureable $X$ we use the decomposition $X=X^+-X^-$, where $X^+:=\max\{X,0\}\ge 0$ and $X^-:=\max\{-X,0\}\ge0$, to obtain two measureable functions $g^+$ and $g^-$ with $X^+=g^+(Y)$ and $X^-=g^-(Y)$ which yields $X=g(Y)$ for $g:=g^+-g^-$.
\end{proof}
\end{lemma}

\begin{remark}
 Let $(\Omega,\calF,P)$ be a probability space and $Z_1$ and $Z_2$ be two random variables with $Z_1=Z_2$ almost surely. In general, a random variable $X$ which is $\sigma(Z_1)$-measureable is not necessarily measureable with respect to $\sigma(Z_2)$, but since $X=g(Z_1)$ for some measureable $g$ by the factorisation lemma, it holds that $X$ almost surely equals the $\sigma(Z_2)$-measureable variable $g(Z_2)$. Hence for $p\ge 1$ we have $L^p(\Omega,\sigma(Z_1),P)=L^p(\Omega,\sigma(Z_2),P)$ as subspaces of $L^p(\Omega,\calF,P)$ as in Remark~\ref{augmentationremark}. This allows us to define $\E\big[X\big|[Z]\big]:=\E[X|Z]$, where $[Z]$ is an equivalence class of random variables with respect to almost sure equality.
\end{remark}

The following, sometimes called \emph{L\'evy's zero-one law}, is an implication of Doob's well-known martingale convergence theorem, see e.g. \cite{Bogachev2007,Oksendal2003}.

\begin{theorem}[L\'evy's zero-one law]
 Let $(\Omega,\calF,P)$ be a probability space, $(\calF_n)_{n\in\N}$ a filtration, $\calF_\infty:=\sigma(\bigcup_{n\in\N}\calF_n)$ and $X\in L^1(\Omega)$. Then $\lim_{n\to\infty}\E[X|\calF_n]=\E[X|\calF_\infty]$ in $L^1(\Omega)$.
\end{theorem}

Since we focus on the space of $L^2$-functions, we need the following proposition:

\begin{proposition}\label{expectconv}
 Let $(\Omega,\calF,P)$ be a probability space, $(\calF_n)_{n\in\N}$ a filtration and assume $X\in L^2(\Omega)$. If $(\E[X|\calF_n])_{n\in\N}$ is a Cauchy sequence in $L^2(\Omega)$, then for $\calF_\infty:=\sigma(\bigcup_{n\in\N}\calF_n)$ we have $\lim_{n\to\infty}\E[X|\calF_n]=\E[X|\calF_\infty]$ in $L^2(\Omega)$.
\begin{proof}
 Let $Y\in L^2(\Omega)$ be the limit of $(\E[X|\calF_n])_{n\in\N}$. Since we have $\|\cdot\|_{L^1}\le\|\cdot\|_{L^2}$ on $L^2(\Omega)$ by Remark~\ref{Hoelderprob}, together with L\'evy's zero-one law we get
\begin{align*}
 \|Y-\E[X|\calF_\infty]\|_{L^1}
 & \le\inf_{n\in\N}\|Y-\E[X|\calF_n]\|_{L^1}+\|\E[X|\calF_n]-\E[X|\calF_\infty]\|_{L^1}\\
 & \le\inf_{n\in\N}\|Y-\E[X|\calF_n]\|_{L^2}+\|\E[X|\calF_n]-\E[X|\calF_\infty]\|_{L^1}=0,\\
\end{align*}
so $Y=\E[X|\calF_\infty]$ in $L^1(\Omega)$. Since $Y\in L^2(\Omega)$, we also get $Y=\E[X|\calF_\infty]$ in $L^2(\Omega)$.
\end{proof}
\end{proposition}

\begin{remark}
 We note that in the above proposition the assumption for $(\E[X|\calF_n])_{n\in\N}$ to be a Cauchy sequence in $L^2$ is actually redundant, since $\left\|\E[X|\calF_n]\right\|_{L^2}\le\|X\|_{L^2}$ for all $n\in\N$, hence $(\E[X|\calF_n])_{n\in\N}$ is a bounded martingale in $L^2$, and one can show that a martingale which is bounded in $L^p$ for some $p\in(1,\infty)$ already converges in $L^p$, see e.g. \cite{Bogachev2007}.
\end{remark}

\begin{theorem}\label{condexpectfinite}
 Let $n\in\N_0$, $m\in\N$, $f^{(n)}\in\ell_A^2(\calH)^\wotn$ and let $\{\psi_1,\dots,\psi_m\}$ be an orthonormal system in $\ell_A^2(\calH)$. Then for $\calG:=\sigma\big(\la\psi_k,\cdot\ra,k=1,\dots,m\big)$ we have
\begin{equation}\label{condexpectfiniteformula}
 \E\left[\wickpoly{f^{(n)}}\Big|\calG\right]=\wickpoly{\calP_\psi^\otn f^{(n)}},
\end{equation}
where $\calP_\psi:\ell_A^2(\calH)\to\spann\{\psi\}:=\spann\{\psi_1,\dots,\psi_m\}$ is the orthogonal projection, i.e.
\[ \calP_\psi f=\sum_{k=1}^m(f,\psi_k)_A\psi_k\quad\text{for }f\in\ell_A^2(\calH). \]
\begin{proof}
 For $n=0$ we have $\wickpoly{f^{(n)}}=\wickpoly{\calP_\psi^\otn f^{(n)}}\in\R$, hence the statement is clear in that case and we may assume $n\not=0$. We will first prove the assertion for $f^{(n)}=\phi^\otn$ for some $\phi\in s(\calN)$ only and afterwards derive that the property transfers to arbitrary $f^{(n)}\in\ell_A^2(\calH)^\wotn$ by density arguments. So let $\phi\in s(\calN)$. First note that $\wickpoly{\calP_\psi^\otn\phi^\otn}=\wickpoly{(\calP_\psi\phi)^\otn}$ is $\calG$-measureable by Corollary~\ref{wickhermite} and we have $(\psi_k,\phi-\calP_\psi\phi)_A=0$ for $k=1,\dots,m$. If $\phi\in\spann\{\psi\}$, then $\wickpoly{\phi^\otn}$ is $\calG$-measureable and hence
\[ \E\left[\wickpoly{\phi^\otn}\Big|\calG\right]=\wickpoly{\phi^\otn}=\wickpoly{\calP_\psi^\otn\phi^\otn}. \]
Otherwise $\phi\not\in\spann\{\psi\}$, so $\phi-\calP_\psi\phi\not=0$. For $G\in\calG$ there exists a measureable function $g$ with $\indi_G=g(\la\psi_1,\cdot\ra,\dots,\la\psi_m,\cdot\ra)$ by Lemma~\ref{factorisationlemma}. To shorten the notation we write $g(\psi)=g(\la\psi_1,\cdot\ra,\dots,\la\psi_m,\cdot\ra)$. We distinguish two cases: If $\calP_\psi\phi=0$, then by Corollary~\ref{wickhermite} we have
\[ \wickpoly{\phi^\otn}=\|\phi\|_A^nH_n\left(\frac{\la\phi,\cdot\ra}{\|\phi\|_A}\right)=\|\phi\|_A^nH_n\left(\frac{\la\phi-\calP_\psi\phi,\cdot\ra}{\|\phi-\calP_\psi\phi\|_A}\right), \]
which together with Corollary~\ref{orthofubinigeneral}, Equation~\eqref{Hermiteintegraliszero} from page~\pageref{Hermiteintegraliszero} and $n\not=0$ yields
\begin{align*}
 \int_G\wickpoly{\phi^\otn}d\mu_A
 & = \|\phi\|_A^n\int_{s'}g(\psi)H_n\left(\frac{\la\phi-\calP_\psi\phi,\cdot\ra}{\|\phi-\calP_\psi\phi\|_A}\right)d\mu_A\\
 & = \|\phi\|_A^n\int_{s'}g(\psi)d\mu_A\int_{\R}H_n(x)d\mu_1(x)\\
 & = \|\phi\|_A^n\mu_A(G)\delta_{0,n}\\
 & = 0\\
 & = \int_G\wickpoly{\calP_\psi^\otn\phi^\otn}d\mu_A.
\end{align*}
In the other case, $\calP_\psi\phi\not=0$, we recall the assumption $\phi-\calP_\psi\phi\not=0$ to observe
\[ 0<\|\calP_\psi\phi\|_A < \|\calP_\psi\phi\|_A+\|\phi-\calP_\psi\phi\|_A=\|\phi\|_A, \]
hence for $\beta:=\|\calP_\psi\phi\|_A\cdot\|\phi\|_A^{-1}$ we have $\beta\in(0,1)$, thus $\alpha:=\sqrt{1-\beta^2}\in(0,1)$ and it holds $\alpha^2+\beta^2=1$. Corollary~\ref{wickhermite} and Equation~\eqref{Hermitebinomial} yield
\begin{align}
 \wickpoly{\phi^\otn}
 & = \|\phi\|_A^nH_n\left(\frac{\la\phi,\cdot\ra}{\|\phi\|_A}\right)\notag\\
 & = \|\phi\|_A^nH_n\left(\frac{\la\phi-\calP_\psi\phi,\cdot\ra}{\|\phi\|_A}+\frac{\la \calP_\psi\phi,\cdot\ra}{\|\phi\|_A}\right)\notag\\
 & = \|\phi\|_A^nH_n\left(\frac{\alpha\la\phi-\calP_\psi\phi,\cdot\ra}{\alpha\|\phi\|_A}+\frac{\beta\la\calP_\psi\phi,\cdot\ra}{\beta\|\phi\|_A}\right)\notag\\
 & = \|\phi\|_A^n\sum_{k=0}^n\binom{n}{k}\alpha^k\beta^{n-k}\underbrace{H_k\left(\frac{\la\phi-\calP_\psi\phi,\cdot\ra}{\alpha\|f\|_A}\right)H_{n-k}\left(\frac{\la \calP_\psi\phi,\cdot\ra}{\beta\|\phi\|_A}\right)}_{=:I_k(\cdot)}\label{somelabel3}.
\end{align}
For $k\in\{0,\dots,n\}$ by Corollary~\ref{orthofubinigeneral} and Equation~\eqref{Hermiteintegraliszero} we have
\begin{align}
 \int_{s'}g(\psi)I_kd\mu_A
 & = \int_{s'}g(\psi)H_{n-k}\left(\frac{\la\calP_\psi\phi,\cdot\ra}{\beta\|\phi\|_A}\right)d\mu_A\cdot\int_{s'}H_k\left(\frac{\la\phi-\calP_\psi\phi,\cdot\ra}{\alpha\|\phi\|_A}\right)d\mu_A\notag\\
 & = \int_{s'}g(\psi)H_{n-k}\left(\frac{\la\calP_\psi\phi,\cdot\ra}{\|\calP_\psi\phi\|_A}\right)d\mu_A\cdot\delta_{0,k},\label{somelabel4}
\end{align}
where we used $\alpha\|\phi\|_A=\|\phi-\calP_\psi\phi\|_A$. Finally we obtain
\begin{align*}
 \int_G\wickpoly{\phi^\otn}d\mu_A
 & = \int_{s'}g(\psi)\wickpoly{\phi^\otn}d\mu_A\\
 & \stackrel{\eqref{somelabel3}}{=} \|\phi\|_A^n\sum_{k=0}^n\binom{n}{k}\alpha^k\beta^{n-k}\int_{s'}g(\psi)I_kd\mu_A\\
 & \stackrel{\eqref{somelabel4}}{=} \int_{s'}g(\psi)\|\calP_\psi\phi\|_A^nH_n\left(\frac{\la\calP_\psi\phi,\cdot\ra}{\|\calP_\psi\phi\|_A}\right)d\mu_A\\
 & = \int_G\wickpoly{\calP_\psi^\otn\phi^\otn}d\mu_A.
\end{align*}
We have established $\E\left[\wickpoly{\phi^\otn}|\calG\right]=\wickpoly{\calP_\psi^\otn\phi^\otn}$ for $\phi\in s(\calN)$. Now for some $\phi^{(n)}\in s(\calN)^\wotn$ with representation
\[ \phi^{(n)}=\sum_{k=1}^m\alpha_k\phi_k^\otn\quad\text{for some }m\in\N,\,\alpha\in\R^m\text{ and }\phi_1,\dots,\phi_m\in\ell_A^2(\calH) \]
by linearity of the conditional expectation we have
\[ \E\left[\wickpoly{\phi^{(n)}}|\calG\right]=\sum_{k=1}^m\alpha_k\wickpoly{\calP_\psi^\otn\phi_k^\otn}=\wickpoly{\calP_\psi^\otn \phi^{(n)}}. \]
Since $s(\calN)^\wotn$ is dense in $\ell_A^2(\calH)^\wotn$ and both $\calP_\psi^\otn$ and the conditional expectation are continuous, we also have \eqref{condexpectfiniteformula} for $f^{(n)}\in\ell_A^2(\calH)^\wotn$.
\end{proof}
\end{theorem}

\begin{corollary}
 Let $F\in L^2(\mu_A)$ with chaos decomposition $F=\sum_{n=0}^\infty\wickpoly{f^{(n)}}$. For $\{\psi_1,\dots,\psi_m\}$, $\calG$ and $\calP_\psi$ as in the previous theorem we have
\begin{equation}
 \E[F|\calG]=\sum_{n=0}^\infty\wickpoly{\calP_\psi^\otn f^{(n)}}.
\end{equation}
\begin{proof}
 The conditional expectation operator is linear and continuous, hence
\[ \E[F|\calG]=\sum_{n=0}^\infty\E\left[\wickpoly{f^{(n)}}|\calG\right]=\sum_{n=0}^\infty\wickpoly{\calP_\psi^\otn f^{(n)}}. \]
\end{proof}
\end{corollary}

\begin{lemma}\label{samespansamesigmaalgebra}
 Let $(h_i)_{i\in\N}$ and $(h'_j)_{j\in\N}$ be two orthonormal bases in $\calH$ and assume the sets $\{x_1,\dots,x_n\}$ and $\{y_1,\dots,y_m\}$ span the same subspace in $\ell^2(\R)$ for some $n,m\in\N$. Then for the $\sigma$-algebras
\[ \calG_1:=\sigma\big(\la h_i\bullet x_k,\cdot\ra:i\in\N,k=1,\dots,n\big)\quad\text{and}\quad\calG_2:=\sigma\big(\la h'_j\bullet y_l,\cdot\ra:j\in\N,l=1,\dots,m\big) \]
we have that $L^2(s'(\calN),\calG_1,\mu_A)=L^2(s'(\calN),\calG_2,\mu_A)$ as subspaces of $L^2(\mu_A)$, i.e. in the sense of Remark~\ref{augmentationremark}.
\begin{proof}
 By symmetry it suffices to show $L^2(s'(\calN),\calG_1,\mu_A)\subset L^2(s'(\calN),\calG_2,\mu_A)$. For $i\in\N$ and $k\in\{1,\dots,n\}$ we have
\[ h_i\bullet x_k=\left(\sum_{j=1}^\infty(h_i,h'_j)h'_j\right)\bullet\left(\sum_{l=1}^m\alpha_ly_l\right)=\sum_{j=1}^\infty\sum_{l=1}^m\alpha_l(h_i,h'_j)h'_j\bullet y_l\quad\text{for some }\alpha\in\R^m. \]
Hence $\la h_i\bullet x_k,\cdot\ra$ is the limit of $\left(\sum_{j=1}^N\sum_{l=1}^m\alpha_l(h_i,h'_j)\la h'_j\bullet y_l,\cdot\ra\right)_{N\in\N}$ in the closed subspace $L^2(s'(\calN),\calG_2,\mu_A)$ and thus an element in the latter itself.
\end{proof}
\end{lemma}

\begin{definition}
 Let $n\in\N$ and $x_1,\dots,x_n\in\ell^2(\R)$. For $F\in L^2(\mu_A)$ we define
\[ \E[F|x_1,\dots,x_n]:=\E[F|\calG],\quad\text{where }\calG:=\sigma\big(\la h_i\bullet x_k,\cdot\ra:i\in\N,k=1,\dots,n\big) \]
for some orthonormal basis $(h_i)_{i\in\N}$ of $\calH$. This notation makes sense since $\E[F|\calG]$ does not depend on the particular choice of $(h_i)_{i\in\N}$, as proven in the lemma above.
\end{definition}

\begin{corollary}\label{samespansamecondexpect}
 Let $F\in L^2(\mu_A)$ and let $\{x_1,\dots,x_n\}$ and $\{y_1,\dots,y_m\}$ span the same subspace in $\ell^2(\R)$ for some $n,m\in\N$. Then
\[ \E[F|x_1,\dots,x_n]=\E[F|y_1,\dots,y_m]. \]
\end{corollary}

From this point we consider $A\in L(\ell^2(\calH))$ to be induced by some self-adjoint and positive definite operator $A\in L(\ell^2(\R))$.

\begin{remark}
 Let $(f_n)_{n\in\N}$ be a family of measureable functions on some measureable space and let $\calF_m:=\sigma(f_i:i=1,\dots,m)$ and $\calF_\infty:=\sigma(f_n:n\in\N)$. Then $\calF_\infty=\sigma\left(\bigcup_{m\in\N}\calF_m\right)$.
\end{remark}

\begin{theorem}\label{condexpect}
 Let $f\in \ell^2(\calH)$, $n\in\N$ and $x_1,\dots,x_n\in\ell^2(\R)$ be orthonormal with respect to $(\cdot,\cdot)_A$. Then we have
\[ \E[\la f,\cdot\ra|x_1,\dots,x_n]=\sum_{k=1}^n\big\la[f,Ax_k]\bullet x_k,\cdot\big\ra=\la Pf,\cdot\ra, \]
where $P:\ell^2(\R)\to\spann\{x_1,\dots,x_n\}$ is the orthogonal projection with respect to $(\cdot,\cdot)_A$. In particular $\E[\la f,\cdot\ra|x_1,\dots,x_n]=\sum_{k=1}^n\E[\la f,\cdot\ra|x_k]$. Note that we require $f\in\ell^2(\calH)$, since $[f,\cdot]$ is not defined for general $f\in\ell_A^2(\calH)$, as we shall see in Example~\ref{noextension}.
\begin{proof}
 Let $(h_i)_{i\in\N}$ be some orthonormal basis of $\calH$ and note that then $(h_i\bullet x_k)_{i\in\N,k=1,\dots,n}$ is an orthonormal system in $\ell_A^2(\calH)$ by Corollary~\ref{miniregel}. We use Proposition~\ref{expectconv} together with the remark above, Theorem~\ref{condexpectfinite} and Proposition~\ref{rechenregeln} to obtain
\begin{align*}
 \E[\la f,\cdot\ra|x_1,\dots,x_n]
 & = \lim_{N\to\infty}\E[\la f,\cdot\ra|\la h_i\bullet x_k,\cdot\ra,i=1,\dots,N,k=1,\dots,n]\\
 & = \lim_{N\to\infty}\left\la\sum_{k=1}^n\sum_{i=1}^N(f,h_i\bullet x_k)_Ah_i\bullet x_k,\cdot\right\ra\\
 & = \sum_{k=1}^n\left\la\lim_{N\to\infty}\sum_{i=1}^N\big([f,Ax_k],h_i\big)h_i\bullet x_k,\cdot\right\ra\\
 & = \sum_{k=1}^n\big\la[f,Ax_k]\bullet x_k,\cdot\big\ra,
\end{align*}
which proves the first equality. For the orthogonal projection $P:\ell^2(\R)\to\spann\{x_1,\dots,x_n\}$ it clearly holds $Px=\sum_{k=1}^n(x,x_k)_Ax_k$ for $x\in\ell^2(\R)$. The identity $f=\sum_{l=0}^\infty f_l\bullet e_l$, continuity of $[\cdot,\cdot]$ and Corollary~\ref{miniregel} give
\begin{align*}
 \sum_{k=0}^n[f,Ax_k]\bullet x_k
 & = \sum_{k=0}^n\left[\sum_{l=0}^\infty f_l\bullet e_l,Ax_k\right]\bullet x_k\\
 & = \sum_{k=0}^n\sum_{l=0}^\infty[f_l\bullet e_l,Ax_k]\bullet x_k\\
 & = \sum_{k=0}^n\sum_{l=0}^\infty(e_l,x_k)_Af_l\bullet x_k\\
 & = \sum_{l=0}^\infty f_l\bullet\left(\sum_{k=0}^n(e_l,x_k)_Ax_k\right)\\
 & = \sum_{l=0}^\infty f_l\bullet Pe_l\\
 & = Pf,
\end{align*}
where we view $P$ as an operator also defined on $\ell^2(\calH)$ as in Theorem~\ref{matrixextension}. Hence
\[ \E[\la f,\cdot\ra|x_1,\dots,x_n]=\sum_{k=1}^n\big\la[f,Ax_k]\bullet x_k,\cdot\big\ra=\la Pf,\cdot\ra. \]
\end{proof}
\end{theorem}

\end{section}

\begin{section}{Examples and Application}\label{appsect}

\begin{example}\label{noextension}
 In this example we will show that $[\cdot,\cdot]$ does not necessarily possess a continuous extension to $\ell_A^2(\calH)\times\ell^2(\R)$. Let $A\in L(\ell^2(\R))$ be the operator uniquely given by
\[ Ae_n=\frac{1}{n^2}e_n\quad\text{for }n\in\N, \]
which exists since $\|Ax\|\le\|x\|$ for $x\in\spann\{e_1,e_2,\dots\}$. Clearly $A$ is self-adjoint and positive definite. Let $h\in\calH$ with $\|h\|=1$ be arbitrary and define $f_n:=\sum_{k=1}^nh\bullet e_k\in\ell^2(\calH)$ for $n\in\N$. Then $(f_n)_{n\in\N}$ is a Cauchy sequence in $\ell^2(\calH)$ with respect to $(\cdot,\cdot)_A$, since for $n,m\in\N$ it holds
\[ \|f_n-f_m\|_A^2=\left(\sum_{k=m+1}^nh\bullet e_k,\sum_{k=m+1}^nh\bullet\frac{1}{k^2}e_k\right)=\sum_{k=m+1}^n\frac{1}{k^2}. \]
Let $x\in\ell^2(\R)$ be given by $x_k=k^{-1}$ for $k\in\N$. Then
\[ \big\|[f_n,x]\big\|=\left\|\sum_{k=1}^n\frac{1}{k}h\right\|=\sum_{k=1}^n\frac{1}{k}\quad\text{for }n\in\N, \]
so the sequence $\big([f_n,x]\big)_{n\in\N}$ is unbounded and hence does not converge in $\calH$. Thus no continuous extension of $[\cdot,\cdot]$ onto $\ell_A^2(\calH)\times\ell^2(\R)$ exists, since otherwise we would have $\lim_{n\to\infty}[f_n,x]=[f,x]\in\calH$, where $f$ is the limit of $(f_n)_{n\in\N}$ in $\ell_A^2(\calH)$.
\end{example}

\begin{example}
 Let the linear operator $A:\ell^2(\R)\to\ell^2(\R)$ be given by
\[ Ae_1=e_1+\frac{1}{2}e_2,\quad Ae_2=\frac{1}{2}e_1+e_2\quad\text{and}\quad Ae_n=e_n\,\text{for }n\ge 3. \]
With respect to $(e_k)_{k\in\N}$ the matrix representation of $A$ becomes
\[ \big((e_k,e_l)_A\big)_{k,l\in\N}=\left(\begin{array}{c|c}
 \begin{matrix}1&\frac{1}{2}\\[1mm]\frac{1}{2}&1\end{matrix} & 0 \\[3mm]\hline
 0 & \Id 
\end{array}\right), \]
which can easily be seen to be bounded, self-adjoint and positive definite. Since only finitely many off-diagonal entries are distinct from zero, we clearly have $\ell_A^2(\calH)=\ell^2(\calH)$. For $f\in\ell^2(\calH)$ Theorem~\ref{condexpect} yields
\begin{enumerate}
 \item $\E[\la f,\cdot\ra|e_1]=\big\la(f_1+\frac{1}{2}f_2)\bullet e_1,\cdot\big\ra$,
 \item $\E[\la f,\cdot\ra|e_2]=\big\la(\frac{1}{2}f_1+f_2)\bullet e_2,\cdot\big\ra$ and
 \item $\E[\la f,\cdot\ra|e_n]=\big\la f_n\bullet e_n,\cdot\big\ra$ for $n\ge 3$.
\end{enumerate}
However, we cannot directly apply the theorem to compute $\E[\la f,\cdot\ra|e_1,e_2]$, since $e_1$ and $e_2$ are not orthogonal with respect to $(\cdot,\cdot)_A$. By defining
\[ \wh{e_2}:=\frac{e_2-(e_2,e_1)_Ae_1}{\|e_2-(e_2,e_1)_Ae_1\|_A}=\sqrt{\frac{4}{3}}\Big(e_2-\frac{1}{2}e_1\Big) \]
we have $(e_1,\wh{e_2})_A=0$, $\|e_1\|_A=\|\wh{e_2}\|_A=1$ and $\spann\{e_1,e_2\}=\spann\{e_1,\wh{e_2}\}$, hence Corollary~\ref{samespansamecondexpect} makes Theorem~\ref{condexpect} applicable:
\begin{align*}
 \E[\la f,\cdot\ra|e_1,e_2]=\E[\la f,\cdot\ra|e_1,\wh{e_2}]
 & = \big\la[f,Ae_1]\bullet e_1,\cdot\big\ra+\big\la[f,A\wh{e_2}]\bullet\wh{e_2},\cdot\big\ra\\
 & = \left\la\left(f_1+\frac{1}{2}f_2\right)\bullet e_1,\cdot\right\ra+\left\la f_2\bullet\left(e_2-\frac{1}{2}e_1\right),\cdot\right\ra\\
 & = \la f_1\bullet e_1,\cdot\ra+\la f_2\bullet e_2,\cdot\ra.
\end{align*}
\end{example}

\begin{application}
 In \cite{FrankSeibold2009}, the authors were engaged with the partial differential equation for radiative transfer, that is	
\begin{equation}\label{rteq}
 \partial_tI(x,\mu,t)+\mu\partial_x(x,\mu,t)+(\sigma(x)+\kappa(x))I(x,\mu,t)=\frac{\sigma(x)}{2}\int_{-1}^1I(x,\mu',t)d\mu'+q(x,t)
\end{equation}
with $t>0$, $x\in(a,b)$ and $\mu\in[-1,1]$. The following approach was used: Set
\[ I_l(x,t):=\int_{-1}^1I(x,\mu,t)P_l(\mu)d\mu=(I(x,\cdot,t),P_l)_{L^2([-1,1])}\quad\text{for }l=0,1,2,\dots, \]
where $P_l$ are the Legendre Polynomials, which form a complete orthogonal system in the Hilbert space $L^2([-1,1])$ and satisfy $\|P_l\|_{L^2([-1,1])}^2=\frac{2}{2l+1}$. Using the recursion relation for the Legendre Polynomials it was proven that \eqref{rteq} is equivalent to the infinite tridiagonal system of first-order partial differential equations
\begin{equation}\label{rteq2}
 \partial_tI_k+b_{k,k-1}\partial_xI_{k-1}+b_{k,k+1}\partial_xI_{k+1}=-c_kI_k+q_k,\quad k=0,1,2,\dots,
\end{equation}
where
\[ b_{k,l}=\frac{k+1}{2k+1}\delta_{k+1,l}+\frac{k}{2k+1}\delta_{k-1,l},\quad c_k=\begin{cases}\kappa & k=0\\\kappa+\sigma & k>0\end{cases},\quad\text{and }q_k=\begin{cases}2\kappa q & k=0\\0 & k>0\end{cases}. \]
In order to start numerical computations, only the first $N$ equations in \eqref{rteq2} can be considered. The problem is to decide how to replace the dependence on $I_{N+1}$ in the equation for $I_N$. A simple approach would be to truncate the system by setting $I_l=0$ for $l>N$, which is called the $P_N$ closure. The approach focussed in \cite{FrankSeibold2009} was the method of \emph{optimal prediction}: Assume one is aware of some correlation between the moments $I_l$, $l=0,1,2,\dots$ via a correlation matrix $A$. Instead of simply neglecting $I_{N+1}$, the information of $I_0,\dots,I_N$ could be used to compute the mean solution for $I_{N+1}$, given $I_0,\dots,I_N$. The formula derived and used in \cite{FrankSeibold2009} was
\begin{equation}\label{finalformula}
 \E[I|I_C]=\E\left[\left.\begin{pmatrix}I_C\\I_F\end{pmatrix}\right|I_C\right]=\begin{pmatrix}I_C\\A_{FC}A_{CC}^{-1}I_C\end{pmatrix}=\begin{pmatrix}\Id_{CC}&0\\A_{FC}A_{CC}^{-1}&0\end{pmatrix}I,
\end{equation}
where $C=\{0,\dots,N\}$, $F=\{N+1,N+2,\dots\}$ and the correlation matrix $A$ and the sequence $I$ are split into corresponding blocks
\[ A=\begin{pmatrix}A_{CC}&A_{CF}\\A_{FC}&A_{FF}\end{pmatrix}\quad\text{and}\quad I=\begin{pmatrix}I_C\\I_F\end{pmatrix}. \]
We are going to justify this notation with our results derived about conditional expectations, of course provided all necessary assumptions are fulfilled. Let $\calN\subset\calH\subset\calN'$ be a Gel'fand triple, which gives rise to a Gel'fand triple $s(\calN)\subset\ell^2(\calH)\subset s'(\calN)$ by Theorem~\ref{seqgelfand}. Let a self-adjoint and positive definite operator $A\in L(\ell^2(\R))$ be given and consider the Gaussian measure $\mu_A$ on $s'(\calN)$ as in Definition~\ref{defmuA}. In consistency with the rest of this thesis, we stick to the agreement $0\not\in\N$, so $C=\{1,\dots,N\}$ and $F=\{N+1,N+2,\dots\}$. We identify $A$ with the infinite matrix $\big((e_k,Ae_l)\big)_{k,l\in\N}$ and note that applying $A$ to a sequence $x\in\ell^2(\R)$ simply becomes usual infinite-dimensional matrix multiplication. Note that $A_{CC}$ is positive definite and thus bijective on $\R^C$ with inverse $A_{CC}^{-1}$. Consider the matrix
\[ P:=\begin{pmatrix}\Id_{CC}&A_{CC}^{-1}A_{CF}\\0&0\end{pmatrix}=\begin{pmatrix}A_{CC}^{-1}&0\\0&0\end{pmatrix}\begin{pmatrix}A_{CC}&A_{CF}\\0&0\end{pmatrix} \]
which defines a linear operator $P:\spann\{e_1,e_2,\dots\}\to\spann\{e_1,\dots,e_N\}$ by infinite matrix multiplication. It can easily be verified that $P$ has a continuous linear extension on $\ell^2(\R)$ with operator norm $\|P\|_{L(\ell^2(\R))}\le\|A_{CC}^{-1}\|_{L(\R^C)}\|A\|_{L(\ell^2(\R))}$. Similarly for
\[ P^T:=\begin{pmatrix}\Id_{CC}&0\\A_{FC}A_{CC}^{-1}&0\end{pmatrix}=\begin{pmatrix}A_{CC}&0\\A_{FC}&0\end{pmatrix}\begin{pmatrix}A_{CC}^{-1}&0\\0&0\end{pmatrix} \]
we have $P^T:\ell^2(\R)\to\ell^2(\R)$ with $\|P^T\|_{L(\ell^2(\R))}\le\|A\|_{L(\ell^2(\R))}\|A_{CC}^{-1}\|_{L(\R^C)}$. Note that the operators $P$ and $P^T$ are adjoint to each other with respect to $(\cdot,\cdot)_{\ell^2(\R)}$. The obvious identity $AP=P^TA$ yields that for all $x,y\in\ell^2(\R)$ we have
\begin{equation}\label{PselfadjforA}
 (x,Py)_A=(x,APy)=(x,P^TAy)=(Px,Ay)=(Px,y)_A.
\end{equation}
For $x\in\ell^2(\R)$ this equation, together with the fact $P^2=P$, yields
\[ \|Px\|_A^2=(Px,Px)_A=(x,P^2x)_A=(x,Px)_A\le\|x\|_A\|Px\|_A, \]
hence $P$ can be extended to a bounded linear operator $P:\ell_A^2(\R)\to\spann\{e_1,\dots,e_N\}$, where $\ell_A^2(\R)$ denotes the completion of $\ell^2(\R)$ with respect to $(\cdot,\cdot)_A$. Then Equation~\eqref{PselfadjforA} extends to hold for $x,y\in\ell_A^2(\calH)$. One easily sees that $P$ is surjective, hence for $x\in\ell_A^2(\R)$ and $y\in\spann\{e_1,\dots,e_N\}$ it holds $y=Py$ and thus
\[ (x-Px,y)_A=(x-Px,Py)_A=(Px-P^2x,y)_A=(Px-Px,y)_A=0, \]
so $P$ is the orthogonal projection from $\ell_A^2(\R)$ onto $\spann\{e_1,\dots,e_N\}$. Let $\wh{e_1},\dots,\wh{e_N}$ be an orthonormal basis of $\spann\{e_1,\dots,e_N\}$ with respect to $(\cdot,\cdot)_A$. For $f\in\ell^2(\calH)$ by Corollary~\ref{samespansamecondexpect} and Theorem~\ref{condexpect} we have
\[ \E[\la f,\cdot\ra|e_1,\dots,e_N]=\E[\la f,\cdot\ra|\wh{e_1},\dots,\wh{e_N}]=\la Pf,\cdot\ra. \]
We are going to justify Equation~\eqref{finalformula} in the sense that
\[ \la P\phi,\omega\ra=\la\phi,P^T\omega\ra\quad\text{for }\phi\in s(\calN),\ \omega\in s'(\calN). \]
To this end, we prove the following three steps:
\begin{enumerate}
 \item $P\phi\in s(\calN)$ for $\phi\in s(\calN)$, so $\la P\phi,\cdot\ra$ is pointwisely defined.
 \item $P^T\omega\in s'(\calN)$ for $\omega\in s'(\calN)$, so the expression $\la\phi,P^T\omega\ra$ makes sense for $\phi\in s(\calN)$.
 \item $\la P\phi,\omega\ra=\la\phi,P^T\omega\ra$ for $\phi\in s(\calN)$ and $\omega\in s'(\calN)$.
\end{enumerate}
For (i) we show that for all $p\in\N_0$ it holds $P\phi\in\ell_p^2(\calN_p)$ for $\phi\in\ell_p^2(\calN_p)$, and that the map $P:\ell_p^2(\calN_p)\to\ell_p^2(\calN_p)$ is bounded. Let $p\in\N_0$ and $\phi\in\ell_p^2(\calN_p)$. For $\psi:=\left(\begin{smallmatrix}A_{CC}&A_{CF}\\0&0\end{smallmatrix}\right)\phi$ it holds
\[ \psi_k=\sum_{i=1}^\infty A_{ki}\phi_i\quad\text{for }k=1,\dots,N\quad\text{and}\quad\psi_k=0\quad\text{for }k>N. \]
Let $\|\cdot\|_p$ and $(\cdot,\cdot)_p$ denote the norm and inner product on $\calN_p$, respectively, and let $(\eta_l)_{l\in\N}$ be an orthonormal basis of $\calN_p$. Then for $k\in\{1,\dots,N\}$ and all $n,m\in\N$ we have
\begin{align*}
 \left\|\sum_{i=m+1}^n A_{ki}\phi_i\right\|_p^2
 & =   \sum_{i,j=m+1}^nA_{ki}A_{kj}(\phi_i,\phi_j)_p\\
 & =   \sum_{l=1}^\infty\sum_{i,j=m+1}^n(e_k,Ae_i)(e_k,Ae_j)(\phi_i,\eta_l)_p(\phi_j,\eta_l)_p\\
 & =   \sum_{l=1}^\infty\left(e_k,A\sum_{i=m+1}^n(\phi_i,\eta_l)_pe_i\right)^2\\
 & \le \|A\|^2\sum_{l=1}^\infty\sum_{i=m+1}^n(\phi_i,\eta_l)_p^2\\
 & =   \|A\|^2\sum_{i=m+1}^n\|\phi_i\|_p^2\\
 & \le \|A\|^2\sum_{i=m+1}^ni^{2p}\|\phi_i\|_p^2,
\end{align*}
thus $\left(\sum_{i=1}^n A_{ki}\phi_i\right)_{n\in\N}$ is a Cauchy sequence in the complete space $\calN_p$ with limit $\psi_k\in\calN_p$. We have established that $\psi$ is a finite sequence in $\calN_p$, hence an element of $\ell_p^2(\calN_p)$. Since the matrix $\left(\begin{smallmatrix}A_{CC}^{-1}&0\\0&0\end{smallmatrix}\right)$ only has finitely many non-zero entries, it also defines a bounded linear operator on $\ell_p^2(\calN_p)$. Thus $P\phi=\left(\begin{smallmatrix}A_{CC}^{-1}&0\\0&0\end{smallmatrix}\right)\psi\in\ell_p^2(\calN_p)$ with $\|P\phi\|_{\ell_p^2(\calN_p)}\le K\|\phi\|_{\ell_p^2(\calN_p)}$ for some constant $K$. This yields that $\phi\in s(\calN)$ implies $P\phi\in s(\calN)$, since $s(\calN)=\bigcap_{p\in\N_0}\ell_p^2(\calN_p)$. Furthermore the map $P:s(\calN)\to s(\calN)$ is continuous. Then, for $\phi\in s(\calN)$, the conditional expectation $\E[\la\phi,\cdot\ra|e_1,\dots,e_N]$ is even pointwisely defined by
\[ \E[\la\phi,\cdot\ra|e_1,\dots,e_N](\omega)=\la P\phi,\omega\ra\quad\text{for }\omega\in s'(\calN). \]
Similarly, for (ii) we show that for all $p\in\N_0$ it holds $P^T\omega\in\ell_{-p}^2(\calN_{-p})$ for $\omega\in\ell_{-p}^2(\calN_{-p})$. For $p=0$ this has already been established, so let $p\ge 1$ and $\omega\in s'(\calN)$. Again, since $\left(\begin{smallmatrix}A_{CC}^{-1}&0\\0&0\end{smallmatrix}\right)$ only has finitely many non-zero entries, it defines a bounded linear operator on $\ell_{-p}^2(\calN_{-p})$, so $\omega':=\left(\begin{smallmatrix}A_{CC}^{-1}&0\\0&0\end{smallmatrix}\right)\omega\in\ell_{-p}^2(\calN_{-p})$. Then $\omega'':=\left(\begin{smallmatrix}A_{CC}&0\\A_{FC}&0\end{smallmatrix}\right)\omega'=P^T\omega$ is a sequence in $\calN_{-p}$, since for $k\in\N$ it holds
\[ \omega''_k=\sum_{i=1}^NA_{ki}\omega'_i\in\calN_{-p}. \]
We now check that $\omega''$ is a sequence in $\ell_{-p}^2(\calN_{-p})$. Let $\|\cdot\|_{-p}$ and $(\cdot,\cdot)_{-p}$ denote the norm and inner product on $\calN_{-p}$, respectively, and let $(\gamma_l)_{l\in\N}$ be an orthonormal basis of $\calN_{-p}$. We have the norm estimate
\begin{align*}
 \|\omega''_k\|_{-p}^2
 & =   \sum_{i,j=1}^NA_{ki}A_{kj}(\omega'_i,\omega'_j)_{-p}\\
 & =   \sum_{l=1}^\infty\sum_{i,j=1}^N(e_k,Ae_i)(e_k,Ae_j)(\omega'_i,\gamma_l)_{-p}(\omega'_j,\gamma_l)_{-p}\\
 & =   \sum_{l=1}^\infty\left(e_k,A\sum_{i=1}^N(\omega'_i,\gamma_l)_{-p}e_i\right)^2\\
 & \le \|A\|^2\sum_{l=1}^\infty\sum_{i=1}^N(\omega'_i,\gamma_l)_{-p}^2\\
 & =   \|A\|^2\sum_{i=1}^N\|\omega'_i\|_{-p}^2\\
 & \le \|A\|^2N^{2p}\sum_{i=1}^Ni^{-2p}\|\omega'_i\|_{-p}^2\\
 & \le \|A\|^2N^{2p}\|\omega'\|_{\ell_{-p}^2(\calN_{-p})}^2.
\end{align*}
Then by the assumption $p\ge 1$ it holds
\[ \|\omega''\|_{\ell_{-p}^2(\calN_{-p})}^2=\sum_{k=1}^\infty k^{-2p}\|\omega''_k\|_{-p}^2\le\|A\|^2N^{2p}\|\omega'\|_{\ell_{-p}^2(\calN_{-p})}^2\sum_{k=1}^\infty k^{-2p}<\infty, \]
so $P^T\omega=\omega''\in\ell_{-p}^2(\calN_{-p})$ with $\|P^T\omega\|_{\ell_{-p}^2(\calN_{-p})}\le K'\|\omega\|_{\ell_{-p}^2(\calN_{-p})}$ for some constant $K'$. This yields that $\omega\in s'(\calN)$ implies $P^T\omega\in s'(\calN)$, since $s'(\calN)=\bigcup_{p\in\N_0}\ell_{-p}^2(\calN_{-p})$. Furthermore the map $P^T:s'(\calN)\to s'(\calN)$ is continuous. For (iii) we note that $P$ and $P^T$ are adjoint to each other with respect to the usual inner product on $\ell^2(\R)$. Then they are also adjoint to each other as operators on $\ell^2(\calH)$, since for $f,g\in\ell^2(\calH)$ we have
\begin{align*}
(f,Pg)_{\ell^2(\calH)}
 & = \lim_{N\to\infty}\left(\sum_{k=1}^Nf_k\bullet e_k,\sum_{l=1}^Ng_l\bullet Pe_l\right)_{\ell^2(\calH)}\\
 & = \lim_{N\to\infty}\sum_{k,l=1}^N(f_k,g_l)_\calH(e_k,Pe_l)_{\ell^2(\R)}\\
 & = \lim_{N\to\infty}\sum_{k,l=1}^N(f_k,g_l)_\calH(P^Te_k,e_l)_{\ell^2(\R)}\\
 & = \lim_{N\to\infty}\left(\sum_{k=1}^Nf_k\bullet P^Te_k,\sum_{l=1}^Ng_l\bullet e_l\right)_{\ell^2(\calH)}\\
 & = (P^Tf,g)_{\ell^2(\calH)}.
\end{align*}
Since in the chain $s(\calN)\subset\ell^2(\calH)\subset s'(\calN)$ we identified $\ell^2(\calH)$ with its topological dual space, the dual pairing of $\phi\in s(\calN)$ and $\omega\in\ell^2(\calH)\subset s'(\calN)$ is realized as $\la\phi,\omega\ra=(\phi,\omega)_{\ell^2(\calH)}$, so
\[ \la P\phi,\omega\ra=(\phi,P^T\omega)_{\ell^2(\calH)}=\la\phi,P^T\omega\ra. \]
For general $\omega\in s'(\calN)$, there exists $p\in\N_0$ such that $\omega\in\ell_{-p}^2(\calN_{-p})$. Let $(\omega_n)_{n\in\N}$ be a sequence in $\ell^2(\calH)$ approximating $\omega$ with respect to $\|\cdot\|_{\ell_{-p}^2(\calN_{-p})}$. Then
\[ \la P\phi,\omega\ra=\lim_{n\to\infty}\la P\phi,\omega_n\ra=\lim_{n\to\infty}\la\phi,P^T\omega_n\ra=\la\phi,P^T\omega\ra. \]
Now that (i), (ii) and (iii) are proven and thus we have
\[ \E[\la\phi,\cdot\ra|e_1,\dots,e_N](\omega)=\la\phi,P^T\omega\ra\quad\text{for } \phi\in s(\calN),\ \omega\in s'(\calN), \]
we have established \eqref{finalformula} in the weak sense
\[ \E[\omega|\omega_C]=P^T\omega=\begin{pmatrix}\Id_{CC}&0\\A_{FC}A_{CC}^{-1}&0\end{pmatrix}\omega. \]
\end{application}
\end{section}

\end{chapter}

\appendix

\begin{chapter}{Appendix}

\begin{section}{Positive Semidefinite Matrices}

In this section, $(\lambda_{kl})_{k,l=1,\dots,n}$ is assumed to be a Hermitian matrix for some $n\in\N$, i.e. for $k,l=1,\dots,n$ we have $\lambda_{kl}=\overline{\lambda_{lk}}\in\C$.

\begin{definition}
 The matrix $(\lambda_{kl})_{k,l=1,\dots,n}$ is called \emph{positive semidefinite}, if we have
\begin{equation}\label{defposdef}
 \sum_{k,l=1}^n\alpha_k\overline{\alpha_l}\lambda_{kl}\ge 0\quad\text{for all }\alpha\in\C^n.
\end{equation}
\end{definition}

\begin{lemma}\label{realposdeflemma}
 If $\lambda_{kl}\in\R$ for $k,l=1,\dots,n$, then \eqref{defposdef} is equivalent to
\begin{equation}\label{realposdef}
 \sum_{k,l=1}^n\alpha_k\alpha_l\lambda_{kl}\ge 0\quad\text{for all }\alpha\in\R^n.
\end{equation}
\begin{proof}
 Clearly \eqref{defposdef} implies \eqref{realposdef}, so assume \eqref{realposdef} holds and let $\alpha\in\C^n$. For $k=1,\dots,n$ denote $a_k:=\Re(\alpha_k)$ and $b_k:=\Im(\alpha_k)$. Then
\begin{align*}
 \sum_{k,l=1}^n\alpha_k\overline{\alpha_l}\lambda_{kl}
 & = \sum_{k,l=1}^n (a_k+ib_k)(a_l-ib_l)\lambda_{kl}\\
 & = \underbrace{\sum_{k,l=1}^n a_ka_l\lambda_{kl}}_{\ge 0\text{ by }\eqref{realposdef}}+\underbrace{i\sum_{k,l=1}^n b_ka_l\lambda_{kl}-i\sum_{k,l=1}^n a_kb_l\lambda_{kl}}_{=0\text{ since }\lambda_{kl}=\lambda_{lk}}+\underbrace{\sum_{k,l=1}^n b_kb_l\lambda_{kl}}_{\ge 0\text{ by }\eqref{realposdef}}\ge0.\\
\end{align*}
\end{proof}
\end{lemma}

By a theorem in \cite{Schur1911} we have the following:

\begin{theorem}
If $(\nu_{kl})_{k,l=1,\dots,n}$ is another Hermitian matrix and both $(\lambda_{kl})_{k,l=1,\dots,n}$ and $(\nu_{kl})_{k,l=1,\dots,n}$ are positive semidefinite, then so is their pointwise product
$(\lambda_{kl}\nu_{kl})_{k,l=1,\dots,n}$.
\end{theorem}

\begin{corollary}\label{expisposdef}
 If $(\lambda_{kl})_{k,l=1,\dots,n}$ is positive semidefinite, then so is $(\exp(\lambda_{kl}))_{k,l=1,\dots,n}$.
\begin{proof}
 Let $\alpha\in\C^n$. For $m=0$ we have
\[ \sum_{k,l=1}^n\alpha_k\overline{\alpha_l}\lambda_{kl}^m=\sum_{k,l=1}^n\alpha_k\overline{\alpha_l}=\left|\sum_{k=1}^n\alpha_k\right|^2\ge 0 \]
and for $m\ge 1$, an obvious inductive use of Schur's theorem above yields \[ \sum_{k,l=1}^n\alpha_k\overline{\alpha_l}\lambda_{kl}^m\ge 0. \]
Hence
 \[ \sum_{k,l=1}^n\alpha_k\overline{\alpha_l}\exp(\lambda_{kl})=\sum_{k,l=1}^n\alpha_k\overline{\alpha_l}\sum_{m=0}^\infty\frac{1}{m!}\lambda_{kl}^m=\sum_{m=0}^\infty\frac{1}{m!}\sum_{k,l=1}^n\alpha_k\overline{\alpha_l}\lambda_{kl}^m\ge0. \]
\end{proof}
\end{corollary}

\end{section}

\begin{section}{Hermite Polynomials}\label{Hermitepolynomial}

\begin{definition}
 For $n\in\N_0$ define the $n^\text{th}$ Hermite polynomial $H_n\in L^2(\R,\mu_1)$ by
\[ \R\ni x\longmapsto H_n(x):=(-1)^n\exp\left(\frac{1}{2}x^2\right)\frac{\text{d}^n}{\text{d}x^n}\exp\left(-\frac{1}{2}x^2\right)\in\R, \]
where $\mu_1$ is the standard Gaussian measure on $\R$. Then $H_n$ is a polynomial of degree $n$ with
\begin{equation}\label{Hermiteinnerprod}
 (H_n,H_m)_{L^2(\mu_1)}=n!\delta_{n,m}\quad\text{for }n,m\in\N_0.
\end{equation}
In particular $\|H_n\|_{L^2(\mu_1)}^2=n!$, and by $H_0\equiv 1$ we have
\begin{equation}\label{Hermiteintegraliszero}
 \int_{\R}H_n(x)d\mu_1(x)=(H_n,H_0)_{L^2(\mu_1)}=\delta_{0,n}.
\end{equation}
Since the set of polynomials is dense in $L^2(\R,\mu_1)$, the Hermite polynomials form a complete orthogonal system in $L^2(\R,\mu_1)$. All these results can be found in \cite{Bogachev1998}, where the Hermite polynomials are introduced in a slightly different way.
\end{definition}

\begin{remark}
 In literature one may also find the definition of the $n^\text{th}$ Hermite polynomial to be
\[ \R\ni x\longmapsto \wh{H}_n(x):=(-1)^n\exp(x^2)\frac{\text{d}^n}{\text{d}x^n}\exp(-x^2)\in\R. \]
A sum representation for these can be found in \cite{Obata1994}:
\begin{equation}\label{physicistshermitesum}
 \wh{H}_n(x)=\sum_{k=0}^{\lfloor n/2\rfloor}\frac{(-1)^kn!}{k!(n-2k)!}(2x)^{n-2k}\quad\text{for }x\in\R.
\end{equation}
These polynomials, further called \emph{physicists Hermite polynomials}, do \emph{not} form an orthogonal system in $L^2(\R,\mu_1)$, but are orthogonal with respect to the probability measure on $\R$ given by the $\dx$-density
\[ \R\ni x\longmapsto\frac{1}{\sqrt{\pi}}\exp(-x^2)\in\R \]
We can link $H_n$ and $\wh{H}_n$ by the identities
\begin{gather}
\label{Hermiterelation1} H_n(x)=2^{-\frac{n}{2}}\wh{H}_n\left(\frac{x}{\sqrt{2}}\right)\quad\text{and}\\
\label{Hermiterelation2} \wh{H}_n(x)=2^{\frac{n}{2}}H_n\left(\sqrt{2}x\right).
\end{gather}
This also yields a representation similar to Equation~\eqref{physicistshermitesum} for our Hermite polynomials:
\begin{equation}\label{Hermitesumrepresentation}
 H_n(x)=\sum_{k=0}^{\lfloor n/2\rfloor}\frac{(-1)^kn!}{2^kk!(n-2k)!}x^{n-2k}\quad\text{for }x\in\R.
\end{equation}
If one considers the analytical extension to $\C$ of the physicists Hermite polynomials, then for $\alpha,\beta\in\C$ with $\alpha^2+\beta^2=1$ one has an expansion of binomial type
\[ \wh{H}_n(\alpha x+\beta y)=\sum_{k=0}^n\binom{n}{k}\alpha^k\beta^{n-k}\wh{H}_k(x)\wh{H}_{n-k}(y)\quad\text{for }x,y\in\R, \]
see \cite{Obata1994}. By straightforward use of Equations~\eqref{Hermiterelation1} and~\eqref{Hermiterelation2} we obtain the analogous formula
\begin{equation}\label{Hermitebinomial}
H_n(\alpha x+\beta y)=\sum_{k=0}^n\binom{n}{k}\alpha^k\beta^{n-k}H_k(x)H_{n-k}(y)\quad\text{for }x,y\in\R
\end{equation}
for the Hermite polynomials, where we follow the convention $0^0:=1$.
\end{remark}

\end{section}

\begin{section}{Proof of Lemma \ref{wickpolynomial}}\label{wickpolynomialproofsection}
\begin{proof}\label{wickpolynomialproof}
For $n=0,1$ the assertion is clear by definition. Now let $n\ge 2$ and assume the claim has been proven for all natural numbers $0,\dots,n-1$. One computes
\begin{align*}
 \colon\omega^\otn\colon
 & = \omega\wot\colon\omega^{\otimes n-1}\colon-(n-1)\tau_A\wot\colon\omega^{\otimes n-2}\colon\\
 & = \omega\wot\sum_{k=0}^{\lfloor\frac{n-1}{2}\rfloor}\frac{(n-1)!(-1)^k}{2^kk!(n-1-2k)!}\tau_A^{\wot k}\wot\omega^{\otimes n-1-2k}\\
 & \qquad\qquad - (n-1)\tau_A\wot\sum_{k=0}^{\lfloor\frac{n-2}{2}\rfloor}\frac{(n-2)!(-1)^k}{2^kk!(n-2-2k)!}\tau_A^{\wot k}\wot\omega^{\otimes n-2-2k}\\
 & = \sum_{k=0}^{\lfloor\frac{n-1}{2}\rfloor}\frac{(n-1)!(-1)^k(n-2k)}{2^kk!(n-2k)!}\tau_A^{\wot k}\wot\omega^{\otimes n-2k}\\
 & \qquad\qquad - \sum_{k=0}^{\lfloor\frac{n-2}{2}\rfloor}(-1)2\frac{(n-1)!(-1)^{k+1}(k+1)}{2^{k+1}(k+1)!(n-2(k+1))!}\tau_A^{\wot k+1}\wot\omega^{\otimes n-2(k+1)}\\
 & = \sum_{k=0}^{\lfloor\frac{n-1}{2}\rfloor}\frac{n!(-1)^k}{2^kk!(n-2k)!}\tau_A^{\wot k}\wot\omega^{\otimes n-2k} - 2\sum_{k=1}^{\lfloor\frac{n-1}{2}\rfloor}\frac{(n-1)!(-1)^kk}{2^kk!(n-2k)!}\tau_A^{\wot k}\wot\omega^{\otimes n-2k}\\
 & \qquad\qquad + 2\sum_{k=1}^{\lfloor\frac{n}{2}\rfloor}\frac{(n-1)!(-1)^kk}{2^kk!(n-2k)!}\tau_A^{\wot k}\wot\omega^{\otimes n-2k}\\
 & = \underbrace{\sum_{k=0}^{\lfloor\frac{n-1}{2}\rfloor}\frac{n!(-1)^k}{2^kk!(n-2k)!}\tau_A^{\wot k}\wot\omega^{\otimes n-2k}}_{=:A} + \underbrace{2\sum_{k=\lfloor\frac{n-1}{2}\rfloor+1}^{\lfloor\frac{n}{2}\rfloor}\frac{(n-1)!(-1)^kk}{2^kk!(n-2k)!}\tau_A^{\wot k}\wot\omega^{\otimes n-2k}}_{=:B}.
\end{align*}
If $n$ is odd, then $\lfloor\frac{n-1}{2}\rfloor=\lfloor\frac{n}{2}\rfloor$, so $B=0$ and the sum in $A$ runs over $k=0,\dots,\lfloor\frac{n}{2}\rfloor$, which is exactly the claim. If $n$ is even, then $\lfloor\frac{n-1}{2}\rfloor+1=\lfloor\frac{n}{2}\rfloor=\frac{n}{2}$ and hence
\begin{align*}
 A+B
 & = A+2\frac{(n-1)!(-1)^{(n/2)}(n/2)}{2^{(n/2)}(n/2)!(n-2(n/2))!}\tau_A^{\wot (n/2)}\wot\omega^{\otimes n-2(n/2)}\\
 & = A+\frac{n!(-1)^{(n/2)}}{2^{(n/2)}(n/2)!(n-2(n/2))!}\tau_A^{\wot (n/2)}\wot\omega^{\otimes n-2(n/2)}\\
 & =\sum_{k=0}^{\lfloor\frac{n}{2}\rfloor}\frac{n!(-1)^k}{2^kk!(n-2k)!}\tau_A^{\wot k}\wot\omega^{\otimes n-2k}.
\end{align*}
We will spare the reader with the proof of the second equation claimed in Lemma~\ref{wickpolynomial}, as it works similarly.
\end{proof}
\end{section}

\end{chapter}

\end{document}